\RequirePackage{ifpdf}
\ifpdf % We are running pdfTeX in pdf mode
\documentclass[pdftex]{sigma}
\else
\documentclass{sigma}
\fi

\usepackage[mathscr]{eucal} %  use \EuScript (\mathcal unchanged)

\numberwithin{equation}{section}

\newtheorem{Theorem}{Theorem}[section]
\newtheorem{Proposition}[Theorem]{Proposition}
\newtheorem{Lemma}[Theorem]{Lemma}

\theoremstyle{definition}
\newtheorem{Definition}[Theorem]{Definition}
\newtheorem{Example}[Theorem]{Example}
\newtheorem{Remark}[Theorem]{Remark}

\begin{document}

\allowdisplaybreaks

\renewcommand{\PaperNumber}{085}

\FirstPageHeading

\ShortArticleName{Dilogarithm Identities for Sine-Gordon
and Reduced Sine-Gordon Y-Systems}

\ArticleName{Dilogarithm Identities for Sine-Gordon\\
and Reduced Sine-Gordon Y-Systems}

\Author{Tomoki NAKANISHI~$^\dag$ and Roberto TATEO~$^\ddag$}

\AuthorNameForHeading{T.~Nakanishi and R.~Tateo}

\Address{$^\dag$~Graduate School of Mathematics, Nagoya University, Nagoya, 464-8604, Japan}
\EmailD{\href{mailto:nakanisi@math.nagoya-u.ac.jp}{nakanisi@math.nagoya-u.ac.jp}}

\Address{$^\ddag$~Dipartimento di Fisica Teorica and INFN,  Universit\`a di Torino,\\
\hphantom{$^\ddag$}~Via P.~Giuria 1, 10125 Torino, Italy}
\EmailD{\href{mailto:tateo@to.infn.it}{tateo@to.infn.it}}

\ArticleDates{Received May 29, 2010, in f\/inal form October 16, 2010;  Published online October 19, 2010}

\Abstract{We study the family of Y-systems and T-systems associated with the
sine-Gordon models and the reduced sine-Gordon models
for the parameter of continued fractions with two terms.
We formulate these systems by cluster algebras,
which turn out to be of f\/inite type,
and prove their periodicities and the associated dilogarithm
identities which have been conjectured earlier.
In particular, this
provides new examples of periodicities of seeds.}

\Keywords{cluster algebras; quantum groups; integrable models}

\Classification{13F60; 17B37}

\vspace{-3mm}

\section{Introduction}

The Y-systems and T-systems appeared in the study of
two-dimensional
integrable $S$-matrix models and integrable lattice models
in 90's.
They play central roles to connect
 these integrable models with conformal f\/ield theories
\cite{Belavin84} through the method
called the {\em thermodynamic Bethe ansatz\/} (TBA)
(e.g., \cite{Kirillov86,Bazhanov90,
Zamolodchikov90,Zamolodchikov91,Klassen90,Klumper92,
Kuniba92,Ravanini93,Kuniba94a,Tateo95}).

Since the introduction of
the {\em cluster algebras}
by Fomin and Zelevinsky \cite{Fomin02, Fomin03a},
it has been gradually recognized that
cluster algebras provide a suitable
framework to study the mathema\-ti\-cal properties
of Y and T-systems.
As a fruitful outcome, the long standing conjectures
of the periodicity of Y-systems \cite{Zamolodchikov91,
Ravanini93,Kuniba94a}
 (as well as
the periodicity of  T-systems)
and the dilogarithm identities \cite{Kirillov90,Bazhanov90,
Kuniba93a,Gliozzi95}
 for the class of Y and T-systems
associated with the quantum af\/f\/ine algebras
have been
proved
partly
by~\cite{Fomin03b, Chapoton05},
and in full generality
by~\cite{Keller08,Inoue10c,
Nakanishi09,Keller10,Inoue10a,Inoue10b},
recently.

In the above proof, the periodicities of the
Y and T-systems  are reformulated as
the {\em periodici\-ties of seeds\/} in the corresponding cluster algebras.
It was shown in  \cite{Inoue10a} that
 the periodicity of seeds in
a cluster algebra associated with a skew symmetric matrix
reduces
to the periodicity of the
 {\em tropical coefficients} (the `principal coef\/f\/icients'
in \cite{Fomin07}),
which is much simpler than the original problem.
In the above examples,
the periodicities at the level of the tropical coef\/f\/icients
are realized
as combinations of the Coxeter transformations
of the $A$-$D$-$E$ root systems and their variations.
{}From this point of view,
one can regard  these periodicities as a natural extension
of the  periodicities of
the Coxeter mutation sequences
in the cluster algebras of f\/inite type
studied by \cite{Fomin03a,Fomin03b,Fomin07}.
Then, it is natural to ask the following question:
{\em ``Are there any other periodicities of seeds?
And, if there are, can we classify them?''}

The purpose of this paper is to present
 a new class of
(inf\/initely many)  periodicities of seeds in cluster algebras.
These cluster algebras correspond to the
Y-systems introduced
in \cite{Tateo95},
where their periodicities
and the associated dilogarithm identities
were also conjectured.
More precisely, these Y-systems consist of two classes.
The f\/irst class are called the {\em sine-Gordon $($SG$)$ Y-systems\/}
and they are
associated with the TBA equation for the sine-Gordon model~\cite{Zamolodchikov79}.
The second class are called the {\em reduced sine-Gordon $($RSG$)$
Y-systems\/}
and, as the name suggests,
they are  associated with the TBA equation
for a certain reduction of the SG model~\cite{Smirnov90,
Bernard90}.
To these cluster algebras, one can apply the method of~\cite{Inoue10a}, and prove their periodicities
and the associated dilogarithm identities.
This is the main result of the paper.
To be more precise, we concentrate on the case
where the `coupling constant'~$\xi$ is a
continued fraction {\em with two terms}, for simplicity.

The result also suggests us a vague perspective
to the second question above.
Namely, the classif\/ication of the periodicities of seeds
may be comparable with the classif\/ication
of the integrable deformations
of rational conformal f\/ield theories.
In particular, those periodicities
we have already known should be just a tip of iceberg.

The content of the paper is as follows.
In Section~\ref{sec:main} we introduce
the Y-systems and T-systems,
then summarize the results of their periodicities
and the associated dilogarithm identities
in both the SG and the RSG cases.
Here, we concentrate on the special case
of the choice of the coupling constant~$\xi$ in~\eqref{eq:case1}.
In Sections~\ref{sec:clusterSG} and~\ref{sec:proofSG}
we give a proof in the SG case.
In Sections~\ref{sec:clusterRSG} and~\ref{sec:proofRSG}
we give a proof in the RSG case.
In Section~\ref{sec:extension}
we extend the result to a little more
general case of~$\xi$ in~\eqref{eq:case2},
i.e., a general continued fraction with two terms.

In the derivation of the main result,
the properties of the {\em tropical Y-systems\/} in
Propositions~\ref{prop:lev2} and~\ref{prop:lev2r}
are crucial, and the proofs are provided
in detail.
The proofs of the rest are rather formal and
repetitions of the formerly studied cases
 \cite{Keller08,Inoue10c,
Nakanishi09,Keller10,Inoue10a,Inoue10b}.
So, instead of repeating similar proofs, we only provide
examples which typically representing
situations.

{\bf Note added.} After the submission of  the paper,
the anonymous referee pointed out us that
the cluster algebras concerned with the SG and
RSG models studied in this paper
 turn out to be the cluster algebras of type $D$ and $A$,
respectively.
Since this is a very important fact, we include it
in a new subsection (Section \ref{subsec:mutseq}).
This means,
somewhat on the contrary to our initial intention and
expectation,
{\em the periodicities of the mutation sequences
studied in this paper are
actually the ones inside the cluster algebras of finite type.}
In particular, the periodicity property itself
is an automatic consequence of the f\/initeness of the number of seeds.
However, the determination of their
precise periods is still  new in the literature;
furthermore, it does not crucially change our perspective of the
richness of the periodicity phenomena of seeds, which should
be uncovered in the future study.

\section{Main results}
\label{sec:main}

\subsection{Results for sine-Gordon Y-systems}

\begin{figure}[t]
\centerline{\includegraphics{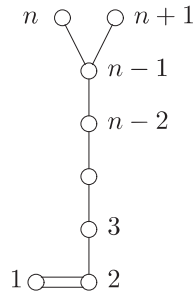}}
\caption{The diagram $X_n$.}
\label{fig:X}
\end{figure}

With an integer $n \geq 4$, we associate a
diagram $X_n$ in Fig.~\ref{fig:X}.
%For later use, we also attach the properties
%$\circ/\bullet$ and $+/-$ to each vertex of $X_n$.
The diagram $X_n$ should {\em not\/} be regarded as an
ordinary Dynkin diagram,
since the horizontality and verticality of
segments also carry information.
It appeared in \cite{Tateo95}
and encodes the structure of the Y-systems
which we are going to introduce now.
Let $\mathcal{I}_n=\{1,\dots,n+1\}\times \mathbb{Z}$.

\begin{Definition}
\label{def:SGY}
Fix an integer $n \geq 4$.
The sine-Gordon (SG)
Y-system $\mathbb{Y}_{n}(\mathrm{SG})$
is the following system of relations for
a family of variables
 $\{Y_i(u) \mid (i,u)\in \mathcal{I}_n \}$,
\begin{gather}
Y_1\left(u-n+1\right)
Y_1\left(u+n-1\right)
=
\left(\prod_{j=2}^{n-1}(1+Y_j(u-n+j))(1+Y_j(u+n-j))\right)\nonumber\\
\phantom{Y_1\left(u-n+1\right)
Y_1\left(u+n-1\right)
=}{}
 \times (1+Y_n(u))(1+Y_{n+1}(u)),
\nonumber\\
Y_2(u-1) Y_2(u+1)
=
\frac
{1+Y_1(u)}
{1+Y_3(u)^{-1}}
,\nonumber\\
Y_i(u-1)Y_i(u+1)
=
\frac{1}{\prod\limits_{j:j\sim i} (1+Y_j(u)^{-1})},
\qquad i=3,\dots,n+1,\label{eq:Y1}
\end{gather}
where $j\sim i$ means that $j$ is adjacent to $i$ in $X_n$.
\end{Definition}

In \cite{Tateo95}, a more general
family of Y-systems was associated with
a rational parameter $\xi$,
which is related  the coupling constant $\beta$
 of the SG model by \eqref{eq:xi}.
The system \eqref{eq:Y1} corresponds to
the special case
\begin{gather}
\label{eq:case1}
\xi=\frac{n-1}{n}=
\cfrac{1}{ 1 +
\cfrac{1}{n-1}
},
\end{gather}
namely, $F=2$, $n_1=1$, and $n_2=n$ in the notation
in \cite{Tateo95}.
The variable $u$ here is related to the variable $\theta$
in \cite{Tateo95} by $u = (2n/\pi\sqrt{-1})\theta$.
Later in Section~\ref{subsec:back} we explain more about the background
of~\eqref{eq:Y1}.

\begin{Definition}
\label{def:YC}
Let $\EuScript{Y}_{n}(\mathrm{SG})$
be the semif\/ield (Appendix \ref{sec:groc}(i)) with generators
$Y_i(u)$
 $((i,u)\in \mathcal{I}_{n})$
and relations $\mathbb{Y}_{n}(\mathrm{SG})$.
Let $\EuScript{Y}^{\circ}_{n}(\mathrm{SG})$
be the multiplicative subgroup
of $\EuScript{Y}_{n}(\mathrm{SG})$
generated by
$Y_i(u)$, $1+Y_i(u)$
 $((i,u)\in \mathcal{I}_{n})$.
(Here we use the symbol $+$ instead of $\oplus$
for simplicity.)
\end{Definition}
The f\/irst main result of the paper is
the following two theorems conjectured by
\cite{Tateo95}.

\begin{Theorem}[Periodicity]
\label{thm:Yperiod}
The following relations hold in $\EuScript{Y}^{\circ}_{n}(\mathrm{SG})$.
\begin{enumerate}\itemsep=0pt
\item[$(i)$] Half periodicity: $Y_i(u+4n-2)=Y_{\omega(i)}(u)$, where
$\omega$ is an involution of the set $\{1,\dots,n+1\}$
defined by $\omega(n)=n+1$, $\omega(n+1)=n$,
and $\omega(i)=i$ $(i=1,\dots,n-1)$.
\item[$(ii)$] Full periodicity: $Y_i(u+8n-4)=Y_i(u)$.
\end{enumerate}
\end{Theorem}

In our proof of Theorem \ref{thm:Yperiod}
we have  a natural interpretation
of the  half period
\[
4n-2=h(D_n)+2+h(D_{n-1})+2
\]
in terms of the Coxeter number  $h(D_n)=2n-2$ of type $D_n$.

Let $L(x)$ be the {\em Rogers dilogarithm function\/}
 \cite{Lewin81}
\begin{gather*}
%\label{eq:L0}
L(x)=-\frac{1}{2}\int_{0}^x
\left\{ \frac{\log(1-y)}{y}+
\frac{\log y}{1-y}
\right\} dy,
\qquad 0\leq x\leq 1.
\end{gather*}
It satisf\/ies the following relation
\begin{gather}
\label{eq:euler}
L(x)+L(1-x)=\frac{\pi^2}{6},
\qquad 0\leq x\leq 1.
\end{gather}

\begin{Theorem}[Functional dilogarithm identities]
\label{thm:dilog}
Suppose that a family of positive real numbers
$\{ Y_i(u) \mid (i,u)\in \mathcal{I}_n\}$ satisfies
$\mathbb{Y}_n(\mathrm{SG})$.
Then, we have the identities
\begin{gather}
\label{eq:DI}
\frac{6}{\pi^2}
\sum_{
\genfrac{}{}{0pt}{1}
{
(i,u)\in \mathcal{I}_{n}
}
{
0\leq u < 8n-4
}
}
L\left(
\frac{Y_i(u)}{1+Y_i(u)}
\right)
 =8(2n-1),\\
\label{eq:DI'}
\frac{6}{\pi^2}
\sum_{
\genfrac{}{}{0pt}{1}
{
(i,u)\in \mathcal{I}_{n}
}
{
0\leq u < 8n-4
}
}
L\left(
\frac{1}{1+Y_i(u)}
\right)
 =4(n-1)(2n-1).
\end{gather}
\end{Theorem}
Two identities
\eqref{eq:DI} and \eqref{eq:DI'}
are equivalent to each other due to \eqref{eq:euler}.

Using this opportunity, we introduce
another system of relations
 accompanying  $\mathbb{Y}_{n}(\mathrm{SG})$.
\begin{Definition}
Fix an integer $n \geq 4$.
The sine-Gordon (SG)
T-system $\mathbb{T}_{n}(\mathrm{SG})$
is the following system of relations for
a family of variables
 $\{T_i(u) \mid (i,u)\in \mathcal{I}_n \}$,
\begin{gather}
T_1\left(u-n+1\right)
T_1\left(u+n-1\right)
=
T_2(u)+1,\nonumber\\
T_2(u-1) T_2(u+1)
=
T_1(u-n+2)T_1(u+n-2)+T_3(u),
\nonumber\\
T_i(u-1)T_i(u+1)
=
T_1(u-n+i)T_1(u+n-i)
+ \prod_{j:j\sim i} T_j(u),
\qquad i=3,\dots,n-1,
\nonumber\\
T_{n}(u-1)T_{n}(u+1)
=
T_1(u) + T_{n-1}(u),
\nonumber\\
T_{n+1}(u-1)T_{n+1}(u+1)
=
T_1(u) + T_{n-1}(u),\label{eq:T1}
\end{gather}
where $j\sim i$ means that $j$ is adjacent to $i$ in $X_n$.
\end{Definition}

There are two connections between
$\mathbb{Y}_{n}(\mathrm{SG})$ and $\mathbb{T}_{n}(\mathrm{SG})$.
The f\/irst connection is a formal one.
Set
\begin{gather}
\label{eq:di}
d_1=n-1, \qquad d_i=1,\qquad i=2,\dots,n-1,
\end{gather}
and let us  write \eqref{eq:Y1}
in a unif\/ied manner as
\begin{gather}
\label{eq:Yuni}
Y_i\left(u-d_i\right)
Y_i\left(u+d_i\right)
=
\frac{
{\prod\limits_{(j,v)\in \mathcal{I}_{n}}
(1+Y_j(v))^{G_+(j,v;i,u)}
}
}
{\prod\limits_{(j,v)\in \mathcal{I}_{n}}
(1+Y_j(v)^{-1})^{G_-(j,v;i,u)}
}.
\end{gather}
Then, \eqref{eq:T1} is written as
\begin{gather}
\label{eq:Tuni}
T_i\left(u-d_i\right)
T_i\left(u+d_i\right)
=
\prod_{(j,v)\in \mathcal{I}_{n}}
T_j(v)^{G_+(i,u;j,v)}
+
\prod_{(j,v)\in \mathcal{I}_{n}}
T_j(v)^{G_-(i,u;j,v)}.
\end{gather}
Note that we took the `transpositions'
of $G_+$ and $G_-$ in \eqref{eq:Tuni}.
The second connection is an algebraic one.
Suppose that $\{ T_i(u) \mid (i,u)
\in \mathcal{I}_n\}$ satisf\/ies
the T-system $\mathbb{T}_n(\mathrm{SG})$.
Set
\[
Y_i(u)
=
\frac{\prod\limits_{(j,v)\in \mathcal{I}_{n}}
T_j(v)^{G_+(i,u;j,v)}
}
{\prod\limits_{(j,v)\in \mathcal{I}_{n}}
T_j(v)^{G_-(i,u;j,v)}}.
\]
Then, $\{ Y_i(u) \mid (i,u)
\in \mathcal{I}_n\}$ satisf\/ies
the Y-system $\mathbb{Y}_n(\mathrm{SG})$.
One may check the claim directly at this moment
using $\mathbb{T}_n(\mathrm{SG})$
(with some ef\/fort).
Alternatively and better, one can automatically obtain  it from
 \cite[Proposition 3.9]{Fomin07}
once we formulate these systems by a cluster algebra
in the next section.

\begin{Definition}
\label{defn:TC}
Let $\EuScript{T}_{n}(\mathrm{SG})$
be the commutative ring over $\mathbb{Z}$
with identity element,  with gene\-ra\-tors
$T_i(u)^{\pm 1}$ $((i,u)\in \mathcal{I}_{n})$
and relations $\mathbb{T}_{n}(\mathrm{SG})$
together with $T_i(u)T_i(u)^{-1}=1$.
Let $\EuScript{T}^{\circ}_{n}(\mathrm{SG})$
be the subring of $\EuScript{T}_{n}(\mathrm{SG})$
generated by
$T_i(u)$ $((i,u)\in \mathcal{I}_{n})$.
\end{Definition}

The following theorem can be proved simultaneously
with Theorem \ref{thm:Yperiod}.

\begin{Theorem}[Periodicity]
\label{thm:Tperiod}
The following relations hold in $\EuScript{T}^{\circ}_{n}(\mathrm{SG})$.
\begin{enumerate}\itemsep=0pt
\item[$(i)$] Half periodicity: $T_i(u+4n-2)=T_{\omega(i)}(u)$, where
$\omega$ is the  one in Theorem~{\rm \ref{thm:Yperiod}}.
\item[$(ii)$] Full periodicity: $T_i(u+8n-4)=T_i(u)$.
\end{enumerate}
\end{Theorem}

\begin{Remark}
Actually,
 $\mathbb{Y}_n(\mathrm{SG})$ and  $\mathbb{T}_n(\mathrm{SG})$
are also considered for $n=3$,
and they coincide with the Y and T-systems of
type $B_2$ with level $2$ in~\cite{Kuniba94a}.
Theorems~\ref{thm:Yperiod}, \ref{thm:dilog},
and~\ref{thm:Tperiod} remain valid for $n=3$
due to~\cite{Inoue10a}.
\end{Remark}

All the results in this subsection will be
extended to a more general case~\eqref{eq:case2}
later in Section~\ref{sec:extension}.

\subsection{Results for reduced sine-Gordon Y-systems}

The SG Y-system in the previous subsection admits
a reduction called the reduced SG Y-system.
It is  obtained from $\eqref{eq:Y1}$ by
formally setting $Y_{n}(u),Y_{n+1}(u)\rightarrow -1$
and $Y_{n-1}(u)\rightarrow \infty$ \cite{Tateo95}.
Let $\tilde{\mathcal{I}}_n=\{1,\dots,n-2\}\times \mathbb{Z}$.

\begin{Definition}
\label{def:RSGY}
Fix an integer $n \geq 4$.
The reduced sine-Gordon (RSG)
Y-system $\mathbb{Y}_{n}(\mathrm{RSG})$
is the following system of relations for
a family of variables
 $\{Y_i(u) \mid (i,u)\in \tilde{\mathcal{I}}_n \}$,
\begin{gather}
Y_1\left(u-n+1\right)
Y_1\left(u+n-1\right)
=
\frac{\prod\limits_{j=2}^{n-2}(1+Y_j(u-n+j))(1+Y_j(u+n-j))
}
{
1+Y_{n-2}(u)^{-1}
},
\nonumber\\
Y_2(u-1) Y_2(u+1)
=
\frac{1+Y_1(u)}{1+Y_3(u)^{-1}},\nonumber\\
Y_i(u-1)Y_i(u+1)
=
\frac{1}
{\prod\limits_{j:j\sim i} (1+Y_j(u)^{-1})
},
\qquad i=3,\dots,n-2,\label{eq:Y1r}
\end{gather}
where $j\sim i$ means that $j \leq n-2$ is adjacent to $i$ in $X_n$.
For $n=4$, the second relation is replaced with
\[
Y_2(u-1) Y_2(u+1)
= 1+Y_1(u).
\]
\end{Definition}

\begin{Definition}
\label{def:YCr}
Let $\EuScript{Y}_{n}(\mathrm{RSG})$
be the semif\/ield with generators
$Y_i(u)$
 $((i,u)\in \tilde{\mathcal{I}}_{n})$
and relations $\mathbb{Y}_{n}(\mathrm{RSG})$.
Let $\EuScript{Y}^{\circ}_{n}(\mathrm{RSG})$
be the multiplicative subgroup
of $\EuScript{Y}_{n}(\mathrm{RSG})$
generated by
$Y_i(u)$, $1+Y_i(u)$
 $((i,u)\in \tilde{\mathcal{I}}_{n})$.
\end{Definition}

The second main result of the paper is
the following two theorems conjectured by
\cite{Tateo95}.
The f\/irst theorem was already proved by \cite{Gliozzi95}
for the RSG Y-systems associated with
a general rational~$\xi$ by using the explicit solution
in terms of cross-ratio.

\begin{Theorem}[{Periodicity \cite{Gliozzi95}}]
\label{thm:Yperiodr}
The following relations hold in $\EuScript{Y}^{\circ}_{n}(\mathrm{RSG})$.

Periodicity: $Y_i(u+4n-2)=Y_i(u)$.
\end{Theorem}

In our proof of Theorem \ref{thm:Yperiodr}
we have  a natural interpretation
of the  period
\[
4n-2=2\{h(A_{n-3})+2+h(A_{n-4})+2\}
\]
in terms of the Coxeter number  $h(A_n)=n+1$ of type $A_n$.

\begin{Theorem}[Functional dilogarithm identities]
\label{thm:dilogr}
Suppose that a family of positive real numbers
$\{ Y_i(u) \mid (i,u)\in \tilde{\mathcal{I}}_n\}$ satisfies
$\mathbb{Y}_n(\mathrm{RSG})$.
Then, we have the identities
\begin{gather}
\label{eq:DIr}
\frac{6}{\pi^2}
\sum_{
\genfrac{}{}{0pt}{1}
{
(i,u)\in \tilde{\mathcal{I}}_{n}
}
{
0\leq u < 4n-2
}
}
L\left(
\frac{Y_i(u)}{1+Y_i(u)}
\right)
=6(2n-5),\\
\label{eq:DIr'}
\frac{6}{\pi^2}
\sum_{
\genfrac{}{}{0pt}{1}
{
(i,u)\in \tilde{\mathcal{I}}_{n}
}
{
0\leq u < 4n-2
}
}
L\left(
\frac{1}{1+Y_i(u)}
\right)
=2(2n^2-11n+17).
\end{gather}
\end{Theorem}
Two identities
\eqref{eq:DIr} and~\eqref{eq:DIr'}
are equivalent to each other due to~\eqref{eq:euler}.

Again, we introduce
the `T-system' accompanying $\mathbb{Y}_{n}(\mathrm{RSG})$.
\begin{Definition}
Fix an integer $n \geq 4$.
The reduced sine-Gordon (RSG)
T-system $\mathbb{T}_{n}(\mathrm{RSG})$
is the following system of relations for
a family of variables
 $\{T_i(u) \mid (i,u)\in \tilde{\mathcal{I}}_n \}$,
\begin{gather}
T_1\left(u-n+1\right)
T_1\left(u+n-1\right)
=
T_2(u)+1,\nonumber\\
T_2(u-1) T_2(u+1)
=
T_1(u-n+2)T_1(u+n-2)+T_3(u),
\nonumber\\
T_i(u-1)T_i(u+1)
=
T_1(u-n+i)T_1(u+n-i)
+ T_{i-1}(u)T_{i+1}(u), \qquad i=3,\dots,n-3,\nonumber\\
T_{n-2}(u-1)T_{n-2}(u+1)
=
T_1(u-2)T_1(u+2)
+ T_1(u) T_{n-3}(u).\label{eq:T1r}
\end{gather}
For $n=4$, the second  and the fourth relations are replaced with
\[
T_2(u-1) T_2(u+1)
=
T_1(u-2)T_1(u+2)+T_1(u).
\]
\end{Definition}

There are connections between
$\mathbb{Y}_{n}(\mathrm{RSG})$ and $\mathbb{T}_{n}(\mathrm{RSG})$
parallel to the ones
between
$\mathbb{Y}_{n}(\mathrm{SG})$ and $\mathbb{T}_{n}(\mathrm{SG})$.

\begin{Definition}
\label{defn:TCr}
Let $\EuScript{T}_{n}(\mathrm{RSG})$
be the commutative ring over $\mathbb{Z}$
with identity element,  with generators
$T_i(u)^{\pm 1}$ $((i,u)\in \tilde{\mathcal{I}}_{n})$
and relations $\mathbb{T}_{n}(\mathrm{RSG})$
together with $T_i(u)T_i(u)^{-1}=1$.
Let $\EuScript{T}^{\circ}_{n}(\mathrm{RSG})$
be the subring of $\EuScript{T}_{n}(\mathrm{RSG})$
generated by
$T_i(u)$ $((i,u)\in \tilde{\mathcal{I}}_{n})$.
\end{Definition}

The following theorem can be proved simultaneously
with Theorem \ref{thm:Yperiodr}.

\begin{Theorem}[Periodicity]
\label{thm:Tperiodr}
The following relations hold in $\EuScript{T}^{\circ}_{n}(\mathrm{RSG})$.

Periodicity: $T_i(u+4n-2)=T_i(u)$.
\end{Theorem}

All the results in this subsection will be
extended to a more general case~\eqref{eq:case2}
later in Section~\ref{sec:extension}.

\subsection{Background in integrable models}
\label{subsec:back}

To  provide the reader with the `big picture behind the scene',
we brief\/ly review
the origins and the consequences
of the Y-systems
(Def\/initions~\ref{def:SGY} and~\ref{def:RSGY}),
the periodicity (Theorems~\ref{thm:Yperiod} and~\ref{thm:Yperiodr}),
and the dilogarithm identities (Theorems~\ref{thm:dilog}
and~\ref{thm:dilogr})  in the context of integrable models.
Mathematically speaking, the  whole content in this subsection
 is completely independent
of the rest of the paper,
so that the reader can safely skip it.

The study of integrable models of quantum f\/ield theory has a long
history, with two initially distinct lines of development.
One comes through classical statistical mechanics, since any lattice model
can be viewed as a regularized Euclidean quantum f\/ield theory; the
other is the direct study of the models as quantum f\/ield theories,
either in Euclidean or Minkowski space.

On the statistical-mechanical side, the subject can be traced back
to Onsager's solution of
the two-dimensional Ising model~\cite{Onsager44} but, much of the modern approach
to these models, owes most to the work of
Baxter  summarized in his book~\cite{Baxter82}.

Viewed directly as a problem in quantum f\/ield theory, the history  begins with the work on the quantum
sine-Gordon model, probably the most famous example of an integrable quantum f\/ield
theory. The Euclidean action of the   sine-Gordon model  is:
\begin{gather}
\label{SGA}
A_{\rm SG} = \int dx^2 \left(\frac{1}{16 \pi} (\partial_{\nu} \varphi)^2 -  2 \mu \cos(\beta \varphi) \right),
\end{gather}
where $0<\beta^2<1$ is related to $\xi$ in \eqref{eq:case1} as
\begin{gather}
\label{eq:xi}
\xi= \frac{\beta^2}{ 1 - \beta^2},
\end{gather}
and $\mu$ f\/ixes the mass scale:   $m \propto \mu^{1/(2-2 \beta^2)}$, where $m$ is the soliton mass.
For $\xi<1$  the theory comprises also solitons-antisoliton bound states~-- the breathers~-- with  masses
\[
m_n=2 m \sin(\pi n \xi/2),\qquad n=1,2,\ldots<1/\xi.
\]
The review by  Zamolodchikov and   Zamolodchikov~\cite{Zamolodchikov79} covers most
of this early work, which focused particularly on characteristics of the
models when def\/ined in inf\/inite spatial volumes.

The main goal when
studying any theory was the
exact calculation of its $S$-matrix, describing the scattering of
arbitrary numbers of elementary excitations.

All of this work concerned massive quantum f\/ield theories,
with f\/inite correlation lengths. By contrast, the initial interest in
statistical mechanics is often the study of models at criticality, where the
correlation length is inf\/inite.

Links between the two approaches began to be built with the
development of conformal f\/ield theory (CFT)~\cite{Belavin84} which showed how powerful
algebraic techniques could be used to solve massless quantum f\/ield
theories, corresponding to the continuum limits of critical lattice models.

Most relevant to the
current context was the subsequent discovery by Zamolodchikov~\cite{Zamolodchikov89}
that
suitable  perturbations, of these conformal f\/ield
theories could lead to models of precisely the sorts
which had been studied previously as massive integrable quantum f\/ield
theories and exact $S$-matrix models. The Euclidean action of a perturbed CFT is
\begin{gather}
A_{\mu}= A_{\rm CFT}+ \mu \int dx^2 \phi(x),
\label{action}
\end{gather}
where $A_{\rm CFT}$ is the action of the  conformal invariant theory and  $\phi$ is a
spinless   primary   f\/ield with conformal dimensions $\Delta=\bar{\Delta}<1$.
The dimensionful coupling $\mu$  measures the
deviation from the critical point and  introduces an independent mass scale
proportional to $\mu^{ 1/(2-2 \Delta)}$. Compa\-ring~(\ref{SGA}) with~(\ref{action}), we see that the sine-Gordon model corresponds
 to the perturbation of a CFT with  central charge $c=1$~-- a free massless boson~-- by  the operator  $\phi_{\rm SG}=2 \cos(\beta \varphi)$
 with conformal  dimension $\Delta_{\rm SG}=\beta^2$.

Another interesting   family of exact $S$-matrix models   is
obtained from the sine-Gordon model, at rational values of~$\beta^2$,
through a quantum group restriction~\cite{Smirnov90, Bernard90} of the Hilbert space.

Setting $\beta^2=p/q$ with $q>p$ coprime integers, this inf\/inite family of models corresponds to the
minimal conformal f\/ield theories $\mathcal{M}_{p,q}$  perturbed by  the operator $\phi_{1,3}$ with
$\Delta_{\rm RSG}=(2p-q)/q$.

As we will see shortly,   it is not the Virasoro central charge  $c=1-6(p-q)^2/pq$  but the ef\/fective central charge
$c_{\rm ef\/f}=c-24 \Delta_{0}$, with $\Delta_{0}$ the dimension of the f\/ield $\phi_{0}$, which generates the ground state of the  conformal
f\/ield theory on a cylinder.  For the  $\mathcal{M}_{p,q}$ family  of models  $c_{\rm ef\/f}=1-6/pq$.

Links with statistical mechanics grew stronger as f\/inite-size ef\/fects
began to be explored, using techniques such as the thermodynamic Bethe
Ansatz (TBA) \cite{Zamolodchikov90}.
In particular, sets of functional relations,
the  Y-systems, began to emerge \cite{Zamolodchikov91} which very closely
paralleled mathe\-ma\-ti\-cal structures such as the `fusion hierarchies' (T-systems)
found in the developments of Baxter's pioneering work on integrable
lattice models~\cite{Kuniba94a}.

In the TBA  approach, the ground state
energy $E_0(R)$ of a massive  integrable quantum f\/ield theory  conf\/ined
on
an inf\/initely long cylinder of circumference  $R$
is written in terms of dressed single-particle energies $\varepsilon_a(\theta)$ (pseudoenergies)
as
\begin{gather}
E_0(R)=- \frac{\pi c_{\rm ef\/f}(r)}{ 6R} =-\frac{1}{2\pi}\sum_{i=1}^{N}\int_{-\infty}^\infty d\theta\; \nu_i(\theta)
\ln\big(1+e^{-\varepsilon_i(\theta)}\big),
\label{en}
\end{gather}
where $\theta$ is the rapidity and  $r$  is related to the  mass $m_1$
of the lightest excitation in the theory and the circumference  $R$ by  $r=R m_1$.  The pseudoenergies  are
the solutions of a  set
of coupled integral equations known as TBA equations. The latter equations  have the general form
\begin{gather}
\varepsilon_i(\theta) = \nu_i(\theta)  -
\sum_{j=1}^N \int_{-\infty}^{\infty}  d\theta \; \phi_{ij}(\theta-\theta') \ln\big(1+e^{-\varepsilon_j(\theta')}\big).
\label{TBA}
\end{gather}
When all scattering is diagonal, the integral equations of the TBA follow directly from
the mass spectrum $\{ m_i \}$ and the two-body   $S$-matrix elements $S_{ij}$:
\[
 \nu_i(\theta)= R m_i \cosh \theta,\qquad \phi_{ij}(\theta)=  \frac{1}{ 2 \pi
 \sqrt{-1} }\frac{d}{d \theta}\ln S_{ij}(\theta).
\]
In this case, the number of pseudoenergies coincides with the number N of particle types in the
original scattering theory. If the scattering is non-diagonal  the TBA derivation  becomes more complicated but
the f\/inal result  can be  still written  in the form~(\ref{TBA}), with some of the  $\nu_i=0$.
The pseudoenergies with  $\nu_i=0$  correspond to  f\/ictitious particles transporting zero energy and zero momentum.
These new particles are often called `magnons', and can be thought of as constructs introduced
to get the counting of states right. The Y functions, the main subject of this paper,
are related to the pseudoenergies as
$Y_i(\theta)=e^{\varepsilon_i(\theta)}$.

The ultraviolet CFT regime  corresponds to  $m_1 \rightarrow 0$ or, equivalently,
to  $r \rightarrow 0$. In this  special limit
\[
E_0(R) \sim - \frac{\pi c_{\rm ef\/f} }{ 6R},
\]
where $c_{\rm ef\/f}=c_{\rm ef\/f}(0)$ is the  ef\/fective   central  charge.

It is  during the  calculation of $c_{\rm ef\/f}$ from equations  (\ref{en}) and (\ref{TBA}) that  sum-rules for the  Rogers dilogarithm function emerge:
\begin{gather}
c_{\rm ef\/f}= \frac{6}{  \pi^2} \sum_{i=1}^N \left[L\left( \frac{1}{ 1+ Y_i}\right)
- L\left( \frac{1}{ 1+ \Upsilon_i}\right) \right].
\label{cc}
\end{gather}
In (\ref{cc}), the  constants $Y_i$ and $\Upsilon_i$
are the stationary values of  $Y_i(\theta)$ in the limits $r \rightarrow 0$ and
$r \rightarrow \infty$, respectively.

The TBA equations contain also exact  information on the conformal dimension $\Delta$ of the perturbing operator $\phi$.
The key idea is to f\/ind a set of functional relations satisf\/ied by the Y functions: a Y-system.
These relations generally imply a periodicity property for the
pseudo\-energies under a certain imaginary shift in $\theta$:
\[
Y_i\big(\theta +  2 \pi \sqrt{-1} P\big) = Y_i(\theta).
\]
The periodicity phenomenon was f\/irst noticed by Al.B.~Zamolodchikov in~\cite{Zamolodchikov91}
and considerations, based on this periodicity suggest that in the far  ultraviolet
region  $c_{\rm ef\/f}(r)$  will have, apart  for a~possible  irregular anti-bulk term,
an expansion in powers of $r^{2/P}$.
This implies for  $\phi$ either the conformal dimension $\Delta=1-1/P$
and an expansion for $c_{\rm ef\/f}(r)$ with both even and odd powers of~$\mu$,
or  $\Delta=1-1/2P$  and an expansion where  only even powers  of~$\mu$ appear.

Let us conclude this  subsection  by demonstrating
the validity of the argument for
 the SG and the RSG models at the simplest case $n=4$ in~\eqref{eq:case1}. This corresponds to
$\xi=3/4$ and $\beta^2=p/q=3/7$.
At this specif\/ic value of coupling, the sine-Gordon  scattering is non-diagonal with a  single soliton-antisoliton
bound state. Setting  $u= 8 \theta/\pi\sqrt{-1}$,
 the corresponding Y-system~is
\begin{gather}
Y_1(u-3) Y_1(u+3)
=
( 1+ Y_5(u)) ( 1+ Y_4(u)) ( 1+ Y_3(u-1))  \nonumber \\
\phantom{Y_1(u-3) Y_1(u+3)
=}{}
\times
( 1+ Y_3(u+1))( 1+ Y_2(u-2)) ( 1+ Y_2(u+2)), \nonumber\\
Y_2(u-1) Y_2(u+1)
= ( 1+ Y_1(u)) \frac{1}{ 1+Y_3(u)^{-1}},\nonumber\\
Y_3(u-1) Y_3(u+1) =\frac{1}{1+Y_2(u)^{-1}}\frac{1}{ 1+Y_4(u)^{-1}} \frac{1}{1+Y_5(u)^{-1}},   \nonumber\\
Y_4(u-1) Y_4(u+1) = \frac{1}{1+Y_3(u)^{-1}},   \qquad
Y_5(u-1) Y_5(u+1) = \frac{1}{1+Y_3(u)^{-1}}.\label{eq:Yst}
\end{gather}

In (\ref{eq:Yst}), $Y_1$ corresponds to  the  breather, $Y_2$
to the soliton-antisoliton pair, $Y_3$, $Y_4$ and $Y_5$ to magnons.
From the point of view of the sine-Gordon ground state energy the
two  magnonic  nodes  $Y_4$ and $Y_5$ are indistinguishable.
Thus,  it is  the half periodicity property
(Theorem~\ref{thm:Yperiod}),
$14$ in $u$, and
 $P=7/8$ in $\theta$, that is relevant to our purposes. Considering also that  the perturbing
f\/ield $\phi_{\rm SG}$ is intrinsically `self-dual' with only  even powers of $\mu$ contributing
to   the expansion of $c_{\rm ef\/f}(r)$,
the f\/inal  result is   $\Delta=1 - 1/2P= 3/7$. This coincides  with the  conformal dimension of~$\phi_{\rm SG}$.

For the dilogarithm sum-rule~\eqref{cc}
 at $r \rightarrow 0$, instead of solving the stationary version of the Y-system
we can simply use  the result \eqref{eq:DI'} averaged over the period:
\[
\frac{6}{\pi^2 } \sum_{i=1}^{n+1}  L\left( \frac{1}{ 1+ Y_i}\right) = \frac{1}{ 8n-4}  \frac{6}{\pi^2 }  \sum_{
\genfrac{}{}{0pt}{1}
{
(i,u)\in \mathcal{I}_{n}
}
{
0\leq u < 8n-4
}
} L\left( \frac{1}{1+ Y_i(u)} \right) =n-1,
\]
with $n=4$.
Similarly, in the limit $r \rightarrow \infty$ both   $Y_1(u)$ and $Y_2(u)$ tend to inf\/inity and decouple. Hence,  for the second
contribution to the sum-rule
(\ref{cc}) we can consider  the simplif\/ied system
\begin{gather*}
\Upsilon_3(u-1) \Upsilon_3(u+1)  = \frac{1}{1+\Upsilon_4(u)^{-1}} \frac{1}{1+\Upsilon_5(u)^{-1}},  \nonumber \\
\Upsilon_4(u-1) \Upsilon_4(u+1)  = \frac{1}{1+\Upsilon_3(u)^{-1}},  \qquad
\Upsilon_5(u-1) \Upsilon_5(u+1)  = \frac{1}{1+\Upsilon_3(u)^{-1}}.%\label{eq:Yst1}
\end{gather*}
The latter  is a $D_3\equiv A_3$ Y-system.
Due to~\cite{Gliozzi95},  the half-period is  $6$
in the variable $u$,  and we have
\[
\frac{6}{\pi^2 } \sum_{i=3}^{5}  L\left( \frac{1}{1+ \Upsilon_i}\right) = \frac{1}{\pi^2 }
\sum_{
\genfrac{}{}{0pt}{1}
{
(i,u)\in D_{3}
}
{
0\leq u < 6
}
} L\left( \frac{1}{1+ \Upsilon_i(u)} \right) = 2.
\]
The result is  $c_{\rm  ef\/f}= 1$, as expected.

The Y-system for the  corresponding RSG model is
\begin{gather*}
Y_1(u-3) Y_1(u+3)
 =
( 1+ Y_2(u-2)) ( 1+ Y_2(u+2)) \frac{1}{1+Y_2(u)^{-1}},  \nonumber\\
Y_2(u-1) Y_2(u+1)
 = 1+ Y_1(u).%\label{eq:Yst3}
\end{gather*}
We still have $P=7/8$ (Theorem \ref{thm:Yperiodr})
 but now  the perturbing operator is `anti self-dual' with respect
 to the  ground state of the conformal f\/ield theory on  an cylinder.
 This means that both even and odd powers of $\mu$ appear in the expansion of~$c_{\rm ef\/f}(r)$.

Therefore,
$\Delta=1 - 1/P= -1/7$, which indeed matches the conformal dimension of  $\phi_{13}$  in~$\mathcal{M}_{3,7}$.
Finally, from  the result~\eqref{eq:DIr'}
we recover the ef\/fective central charge of $\mathcal{M}_{3,7}$:
\[
c_{\rm ef\/f}=   \frac{3}{7 \pi^2}
\sum_{
\genfrac{}{}{0pt}{1}
{
(i,u)\in \tilde{\mathcal{I}}_{4}
}
{
0\leq u < 14
}
}
L\left( \frac{1}{1+ Y_i(u)} \right) = \frac{5}{7}.
\]

\section{Cluster algebras for SG Y-systems}
\label{sec:clusterSG}

In this section we identify $\mathbb{T}_{n}(\mathrm{SG})$ and
$\mathbb{Y}_{n}(\mathrm{SG})$ as relations for cluster variables
and coef\/f\/icients of the cluster algebra
associated with a certain quiver $Q_{n}(\mathrm{SG})$.
We follow \cite{Fomin07}
 for def\/initions and conventions concerning cluster algebras
with coef\/f\/icients,
which are summarized in Appendix \ref{sec:groc} for the reader's
convenience.

\subsection{Parity decompositions of T and Y-systems}

For $(i,u)\in \mathcal{I}_{n}$,
we set the parity condition $\mathbf{P}_{+}$ by,
for even $n$,
\begin{gather*}
%\label{eq:GPcond1}
\mathbf{P}_{+}: \
\begin{cases}
\mbox{$i+u$ is even},& i=1,\dots,n,\\
\mbox{$n+u$ is even}, & i=n+1,
\end{cases}
\end{gather*}
and, for odd $n$,
\begin{gather*}
%\label{eq:GPcond2}
\mathbf{P}_{+}: \
\begin{cases}
\mbox{$u$ is even}, & i=1,\\
\mbox{$i+u$ is even}, & i=2,\dots,n, \\
\mbox{$n+u$ is even}, & i=n+1.\\
\end{cases}
\end{gather*}
Let $\mathbf{P}_-$ be the negation of $\mathbf{P}_+$.
We write, for example, $(i,u):\mathbf{P}_{+}$ if $(i,u)\in
\mathcal{I}_n$ satisf\/ies
$\mathbf{P}_{+}$.
Let $\mathcal{I}_{n\varepsilon}$ ($\varepsilon=\pm$)
be the set of all $(i,u):\mathbf{P}_{\varepsilon}$.
Def\/ine $\EuScript{T}^{\circ}_{n}(\mathrm{SG})_{\varepsilon}$
($\varepsilon=\pm$)
to be the subring of $\EuScript{T}^{\circ}_{n}(\mathrm{SG})$
generated by
 $T_i(u)$
$((i,u)\in \mathcal{I}_{n\varepsilon})$.
Then, we have
$\EuScript{T}^{\circ}_{n}(\mathrm{SG})_+
\simeq
\EuScript{T}^{\circ}_{n}(\mathrm{SG})_-
$
by $T_i(u)\mapsto T_i(u+1)$ and
\begin{gather*}
%\label{eq:Tfact}
\EuScript{T}^{\circ}_{n}(\mathrm{SG})
\simeq
\EuScript{T}^{\circ}_{n}(\mathrm{SG})_+
\otimes_{\mathbb{Z}}
\EuScript{T}^{\circ}_{n}(\mathrm{SG})_-.
\end{gather*}

For  $(i,u)\in \mathcal{I}_{n}$,
we introduce another parity condition $\mathbf{P}'_{+}$ by
\begin{gather*}
%\label{eq:GP'cond1}
\mathbf{P}'_{+}: \
\begin{cases}
\mbox{$i+u$ is odd},& i=1,\dots,n,\\
\mbox{$n+u$ is odd}, & i=n+1,
\end{cases}
\end{gather*}
We have
\[
(i,u):\ \mathbf{P}'_+ \ \Longleftrightarrow\
  (i,u\pm d_i): \ \mathbf{P}_+,
\]
where $d_i$ is given in \eqref{eq:di}.
Let $\mathbf{P}'_-$ be the negation of $\mathbf{P}'_+$.
Let $\mathcal{I}'_{n\varepsilon}$
($\varepsilon=\pm$)  be
 the set of all $(i,u):\mathbf{P}'_{\varepsilon}$.
Def\/ine $\EuScript{Y}^{\circ}_{n}(\mathrm{SG})_{\varepsilon}$
($\varepsilon=\pm$)
to be the subgroup of $\EuScript{Y}^{\circ}_{n}(\mathrm{SG})$
generated by
$Y_i(u)$, $1+Y_i(u)$
$((i,u)\in \mathcal{I}'_{n\varepsilon})$.
Then, we have
$\EuScript{Y}^{\circ}_{n}(\mathrm{SG})_+
\simeq
\EuScript{Y}^{\circ}_{n}(\mathrm{SG})_-
$
by $Y_i(u)\mapsto Y_i(u+1)$,
$1+Y_i(u)\mapsto 1+Y_i(u+1)$,
 and
\begin{gather*}
%\label{eq:Yfact}
\EuScript{Y}^{\circ}_{n}(\mathrm{SG})
\simeq
\EuScript{Y}^{\circ}_{n}(\mathrm{SG})_+
\times
\EuScript{Y}^{\circ}_{n}(\mathrm{SG})_-.
\end{gather*}

From now on, we mainly treat the $+$ parts,
$\EuScript{T}^{\circ}_{n}(\mathrm{SG})_+$
and
$\EuScript{Y}^{\circ}_{n}(\mathrm{SG})_+$.

\subsection[Quiver $Q_n(\mathrm{SG})$]{Quiver $\boldsymbol{Q_n(\mathrm{SG})}$}
\label{subsect:quiver}

%%%%%%%%%%%%%%%%%%%%%%%%%%%%%
% SG
\begin{figure}
\centerline{\includegraphics{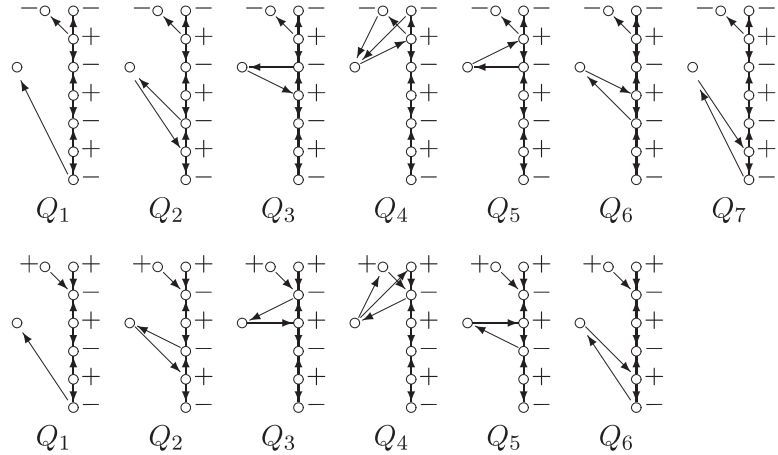}}
\caption{The quiver $Q_{n}(\mathrm{SG})$  for $n=8$
(upper) and for $n=7$ (lower),
where, except for the leftmost vertex of each quiver $Q_i$,
 all the  vertices in the same position
 in  $n-1$ quivers $Q_1,\dots,Q_{n-1}$ are
identif\/ied.}
\label{fig:quiverG}
\end{figure}

Recall that a {\em quiver\/} is an oriented graph,
namely, it consists of the vertices and the arrows
connecting them.
With each $n \geq 4$ we associate
 a quiver $Q_{n}(\mathrm{SG})$ as below.
First, as rather general examples, the cases $n=8$ and $n=7$ are given in
 Fig.~\ref{fig:quiverG},
where, except for the leftmost vertex of each quiver $Q_i$,
 all the  vertices in the same position
 in  $n-1$ quivers $Q_1,\dots,Q_{n-1}$ are
identif\/ied.
For a general $n$,  the quiver $Q_{n}(\mathrm{SG})$ is def\/ined
by naturally extending these examples.
Namely, we consider $n-1$ quivers $Q_1,\dots, Q_{n-1}$,
whose vertices are naturally identif\/ied with
the vertices of the graph $X_n$ in Fig.~\ref{fig:X}.
(The leftmost vertex corresponds to the vertex~$1$ in~$X_n$.)
The arrows are put in $Q_i$ as clearly indicated by the examples
in Fig.~\ref{fig:quiverG}.
Note that the pattern of arrows slightly
depends on the parity of~$n$.
Then, except for the leftmost vertex of each quiver~$Q_i$,
 all the  vertices in the same position
 in  $n-1$ quivers $Q_1,\dots,Q_{n-1}$ are
identif\/ied.
Also we assign the property $+/-$
to  each vertex, except for the leftmost one in each~$Q_i$,
 as in Fig.~\ref{fig:quiverG}.

Let us choose  the index set $\mathbf{I}$
of the vertices of $Q_{n}(\mathrm{SG})$
so that $\mathbf{i}=(i,1)\in \mathbf{I}$ represents
the leftmost vertex in $Q_i$ for $i=1,\dots,n-1$,
 and $\mathbf{i}=(n,i')\in \mathbf{I}$ represents
the vertex $i'=2,\dots,n+1$ in any quiver $Q_i$
under the natural
identif\/ication with $X_n$.
Thus, $i=1,\dots,n$; and $i'=1$ if $i\neq n$
and $i'=2,\dots,n+1$ if $i=n$.

Let $\mathbf{I}_+$
(resp.\ $\mathbf{I}_-$)
denote the set of the vertices $\mathbf{i}\in \mathbf{I}$ with
property  $+$ (resp.\  $-$).
We def\/ine composite mutations (Appendix~\ref{sec:groc}(ii)--(v)),
\begin{gather*}
%\label{eq:mupm2}
\mu_+=\prod_{\mathbf{i}\in\mathbf{I}_+}
\mu_{\mathbf{i}},
\qquad
\mu_-=\prod_{\mathbf{i}\in\mathbf{I}_-}
\mu_{\mathbf{i}}.
\end{gather*}
Note that they do not depend on the order of the product.

For a permutation $\sigma$ of $\{1,\dots,n-1 \}$,
let $\tilde{\sigma}$ be the permutation of
$\mathbf{I}$ such that
$\tilde{\sigma}(i,1)=(\sigma(i),1)$ for $i\neq n$
and $\tilde{\sigma}(n,i')=(n,i')$.
Let $\tilde{\sigma}(Q_{n}(\mathrm{SG}))$
denote the
quiver induced from $Q_{n}(\mathrm{SG})$
 by~$\tilde{\sigma}$.
Namely, if there is an arrow
 $\mathbf{i}\rightarrow
\mathbf{j}$ in $Q_{n}(\mathrm{SG})$,
then, there is an arrow
$\tilde{\sigma}(\mathbf{i})
\rightarrow
\tilde{\sigma}(\mathbf{j})
$
in  $\tilde{\sigma}(Q_{n}(\mathrm{SG}))$.

\begin{Lemma}
\label{lem:GQmut}
Let $Q(0):=Q_{n}(\mathrm{SG})$.
We have the following periodic sequence of mutations of quivers:
\begin{gather}
% 1st
Q(0)
\mathop{\longleftrightarrow}^{{\mu_+} \mu_{(1,1)}}
 Q(1)
\mathop{\longleftrightarrow}^{\mu_-}
Q(2)
\displaystyle
\mathop{\longleftrightarrow}^{\mu_+ \mu_{(2,1)}}
Q(3)
\displaystyle
\mathop{\longleftrightarrow}^{\mu_-}
Q(4)
\nonumber\\
\hphantom{Q(0)
\mathop{\longleftrightarrow}^{{\mu_+} \mu_{(1,1)}}}{}
\mathop{\longleftrightarrow}^{\mu_+\mu_{(3,1)}}
\ \cdots\
\displaystyle
\mathop{\longleftrightarrow}^{\mu_+\mu_{(n-1,1)}}
Q(2n-3)
\displaystyle
\mathop{\longleftrightarrow}^{\mu_-}
Q(2n-2)=Q(0),\label{eq:GB2}
\end{gather}
where the quiver $Q(2p)$ $(p=1,\dots,n-2)$ is given by
\begin{gather}
\label{eq:Q2p}
Q(2p)=
\tilde{\sigma}^p(Q(0)),
\qquad
\sigma=\left(
\begin{matrix}
1&2&\dots &n-2 &n-1\\
2&3& \dots&n-1& 1\\
\end{matrix}
\right).
\end{gather}
\end{Lemma}

\begin{Example}
The mutation sequence
 \eqref{eq:GB2} for $n=6$ is explicitly given
in Fig.~\ref{fig:labelxG1}.
\begin{figure}[t]
\centerline{\includegraphics{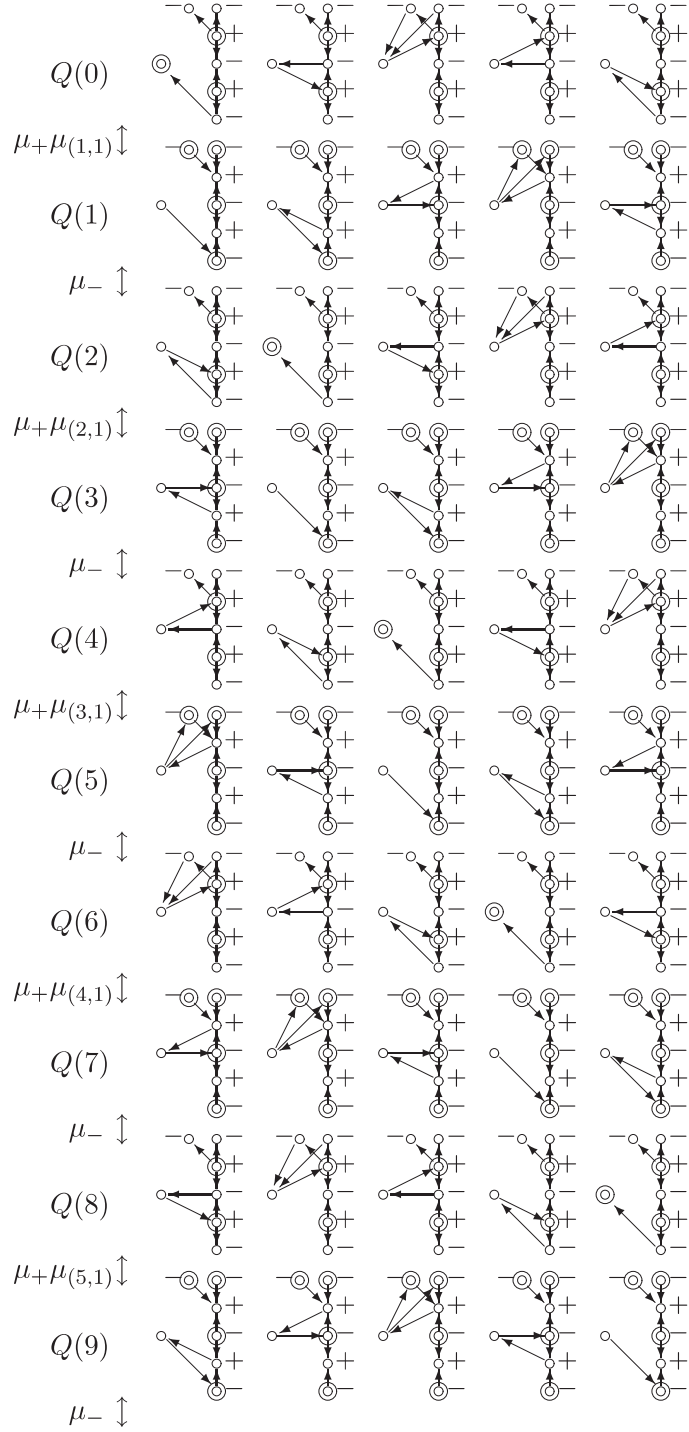}}

\caption{The mutation sequence of the quiver $Q_{n}(\mathrm{SG})$
in~\eqref{eq:GB2} for $n=6$.
The encircled vertices correspond to
the mutation points
$(\mathbf{i},u): \mathbf{p}_+$ in the forward direction.}
\label{fig:labelxG1}
\end{figure}
\end{Example}

\subsection{Embedding maps}
Let
 $B=B_{n}(\mathrm{SG})$
be the skew-symmetric matrix corresponding
to the quiver $Q_{n}(\mathrm{SG})$ (Appendix~\ref{sec:groc}(iii)).
Let $\mathcal{A}(B,x,y)$
be the cluster algebra
 with coef\/f\/icients
in the universal
semif\/ield
$\mathbb{P}_{\mathrm{univ}}(y)$ (Appendix~\ref{sec:groc}(vi)).

In view of Lemma~\ref{lem:GQmut}
we set $x(0)=x$, $y(0)=y$ and def\/ine
clusters $x(u)=(x_{\mathbf{i}}(u))_{\mathbf{i}\in \mathbf{I}}$
 ($u\in \mathbb{Z}$)
 and coef\/f\/icient tuples $y(u)=(y_\mathbf{i}(u))_{\mathbf{i}\in \mathbf{I}}$
 ($u\in \mathbb{Z}$)
by the sequence of mutations
\begin{gather}
% 1st
\cdots
\ \mathop{\longleftrightarrow}^{\mu_-}\
(B(0),x(0),y(0))
\ \mathop{\longleftrightarrow}^{\mu_+
\mu_{(1,1)}}\
(B(1),x(1),y(1))\nonumber\\
\phantom{\cdots}{} \ \mathop{\longleftrightarrow}^{\mu_-}\
\cdots
\ \mathop{\longleftrightarrow}^{\mu_-}\
(B(2n-2),x(2n-2),y(2n-2))
\ \mathop{\longleftrightarrow}^{\mu_+
\mu_{(1,1)}}\
\cdots,\label{eq:Gmutseq}
\end{gather}
where $B(u)$ is the skew-symmetric matrix corresponding to
$Q(u)$.

For  $(\mathbf{i},u)\in
 \mathbf{I}\times \mathbb{Z}$,
we set the parity condition $\mathbf{p}_+$ by
\begin{gather}
\label{eq:GQparity}
\mathbf{p}_+: \
\begin{cases}
 \mathbf{i}\in
\mathbf{I}_+ \sqcup
\{(j+1,1)\},
& u\equiv 2j
\ (j=0,1,\dots,n-2),\\
 \mathbf{i}\in
\mathbf{I}_-,
& \mbox{$u$: odd},
\end{cases}
\end{gather}
where $\equiv$ is modulo $(2n-2)\mathbb{Z}$.
We def\/ine the condition $\mathbf{p}_-$
by $(\mathbf{i},u):\mathbf{p}_- \Longleftrightarrow
(\mathbf{i},u-1):\mathbf{p}_+$.
Plainly speaking,
each $(\mathbf{i},u):\mathbf{p}_+$
(resp.\ $(\mathbf{i},u):\mathbf{p}_-$)
is a mutation point of~\eqref{eq:Gmutseq} in the forward
(resp. backward) direction of~$u$.
See Fig.~\ref{fig:labelxG1}.

\begin{Lemma}
\label{lem:gmap}
Below $\equiv$ means the equivalence modulo $(2n-2)\mathbb{Z}$.

\begin{enumerate}\itemsep=0pt
\item[$(i)$]
The map
\begin{gather*}
g:  \ \mathcal{I}_{n+} \rightarrow
\{ (\mathbf{i},u)\in \mathbf{I}\times \mathbb{Z} \mid
(\mathbf{i},u): \mathbf{p}_+
\}, \\
\phantom{g: \ } \ (i,u-d_i) \mapsto
\begin{cases}
((j+1,1),u), & i= 1; \
u \equiv 2j
\ (j=0,1,\dots,n-2),\\
((n,i),u), & i=2,\dots,n+1
\end{cases}
\end{gather*}
is a bijection.

\item[$(ii)$]
The map
\begin{gather*}
g': \  \mathcal{I}'_{n+} \rightarrow
\{ (\mathbf{i},u)\in \mathbf{I}\times \mathbb{Z} \mid
 (\mathbf{i},u): \mathbf{p}_+
\}, \\
\phantom{g': \ }{} \  (i,u) \mapsto
\begin{cases}
((j+1,1),u), & i= 1;\
u\equiv 2j
\ (j=0,1,\dots,n-2), \\
((n,i),u), &  i=2,\dots,n+1
\end{cases}
\end{gather*}
is a bijection.
\end{enumerate}
\end{Lemma}

Based on Lemma \ref{lem:gmap},
we introduce alternative notations
$\tilde{x}_i(u-d_i):=x_{\mathbf{i}}(u)$
for $(i,u-d_i)\in \mathcal{I}_{n+}$
with $(\mathbf{i},u)=g((i,u-d_i))$
and
$y_i(u):=y_{\mathbf{i}}(u)$
for $(i,u)\in \mathcal{I}'_{n+}$
with $(\mathbf{i},u)=g'((i,u))$,
respectively,
which turn out to be useful.

 Let $\mathcal{A}(B,x)$ be the cluster algebra
with trivial coef\/f\/icients, where $(B,x)$ is
the initial seed
and the coef\/f\/icient semif\/ield
is the {\em trivial semifield\/}
 $\mathbf{1}=\{1\}$ (Appendix~\ref{sec:groc}(i)).
Let $\pi_{\mathbf{1}}:
\mathbb{P}_{\mathrm{univ}}(y)\rightarrow
\mathbf{1}$, $y_{\mathbf{i}}\mapsto 1$ be the projection.
Let $[x_{\mathbf{i}}(u)]_{\mathbf{1}}$
denote the image of $x_{\mathbf{i}}(u)$
 by the algebra homomorphism
$\mathcal{A}(B,x,y)\rightarrow \mathcal{A}(B,x)$
 induced from $\pi_{\mathbf{1}}$.
It is called the {\em trivial evaluation}.

The following lemma follows from the exchange relation
of
cluster variables \eqref{eq:clust}
and the property of the sequence \eqref{eq:GB2}
 observed in Fig.~\ref{fig:labelxG1}.

\begin{Lemma}
\label{lem:Gx2}
Let $G_+$ and $G_-$ be the ones in \eqref{eq:Yuni} and \eqref{eq:Tuni}.
The family $\{\tilde{x}_i(u)
\mid (i,u)\in \mathcal{I}_{n+}\}$
satisfies a system of relations
\begin{gather}
\tilde{x}_i\left(u-d_i\right)
\tilde{x}_i\left(u+d_i\right)
 =
\frac{y_i(u)}{1+y_i(u)}
\prod_{(j,v)\in \mathcal{I}_{n+}}
\tilde{x}_j(v)^{G_+(i,u;j,v)}
\nonumber\\
\phantom{\tilde{x}_i\left(u-d_i\right)
\tilde{x}_i\left(u+d_i\right)
 =}{}
+
\frac{1}{1+y_i(u)}
\prod_{(j,v)\in \mathcal{I}_{n+}}
\tilde{x}_j(v)^{G_-(i,u;j,v)},\label{eq:xy1}
\end{gather}
where $(i,u)\in \mathcal{I}'_{n+}$.
In particular,
the family $\{ [\tilde{x}_i(u)]_{\mathbf{1}}
\,|\, (i,u)\in \mathcal{I}_{n+}\}$
satisfies the T-system $\mathbb{T}_{n}(\mathrm{SG})$
in $\mathcal{A}(B,x)$
by replacing $T_i(u)$ with $[\tilde{x}_i(u)]_{\mathbf{1}}$.
\end{Lemma}

\begin{Example}
Consider the case $n=6$ in Fig.~\ref{fig:labelxG1}.
Let us consider the mutation at the vertex~$(1,1)$ in $Q(0)$,
to which  the variable $\tilde{x}_1(-5)$ is attached.
The next time~$(1,1)$ is mutated is in~$Q(10)$,
where~$\tilde{x}_1(5)$ is attached.
Meanwhile, the only vertex connected to~$(1,1)$ in~$Q(0)$ is~$(6,2)$,
and the  variable attached to $(6,2)$ in $Q(0)$ is equal to
the variable  $\tilde{x}_2(0)$ attached to~$(6,2)$ in~$Q(1)$.
Taking account of the directions of the arrows,
we have the  relation
\begin{align*}
\tilde{x}_1(-5)\tilde{x}_1(5)=
\frac{y_1(0)}{1+y_1(0)} \tilde{x}_2(0)
+ \frac{1}{1+y_1(0)} ,
\end{align*}
which agrees with
\eqref{eq:T1} and
\eqref{eq:xy1}.
Similarly, consider the mutation at the vertex $(6,2)$ in $Q(1)$,
to which  the variable $\tilde{x}_2(0)$ is attached.
The next time $(6,2)$ is mutated is in $Q(3)$,
where $\tilde{x}_2(2)$ is attached.
Meanwhile, the vertices connected to $(6,2)$ in $Q(1)$
are $(1,1)$, $(2,1)$, and $(6,3)$,
and the  variable attached to them are equal to
$\tilde{x}_1(5)$, $\tilde{x}_1(-3)$, and $\tilde{x}_3(1)$, respectively.
Taking account of the directions of the arrows,
we have the relation
\begin{align*}
\tilde{x}_2(0)\tilde{x}_2(2)=
\frac{y_2(1)}{1+y_2(1)} \tilde{x}_1(-3)\tilde{x}_1(5)
+ \frac{1}{1+y_2(1)} \tilde{x}_3(1).
\end{align*}
As the last example,
consider the mutation at the vertex $(6,6)$ in $Q(1)$,
to which  the va\-riab\-le~$\tilde{x}_6(0)$ is attached.
The next time $(6,6)$ is mutated is in $Q(3)$,
where $\tilde{x}_6(2)$ is attached.
Meanwhile, the vertices connected to $(6,6)$ in $Q(1)$
are $(4,1)$ and $(6,5)$,
and the  variable attached to them are equal to
$\tilde{x}_1(1)$ and  $\tilde{x}_5(1)$, respectively.
Taking account of the directions of the arrows,
we have the relation
\[
\tilde{x}_6(0)\tilde{x}_6(2)=
\frac{y_6(1)}{1+y_6(1)} \tilde{x}_1(1)
+ \frac{1}{1+y_6(1)} \tilde{x}_5(1).
\]
The other cases can be checked in similar manners.
\end{Example}

\begin{Definition}
\label{def:T-sub}
The {\em T-subalgebra
$\mathcal{A}_T(B,x)$
of ${\mathcal{A}}(B,x)$
associated with the sequence \eqref{eq:Gmutseq}}
is the subring of
${\mathcal{A}}(B,x)$
generated by
$[x_{\mathbf{i}}(u)]_{\mathbf{1}}$
($(\mathbf{i},u)\in \mathbf{I}\times \mathbb{Z}$),
or equivalently,
generated by
$[\tilde{x}_{i}(u)]_{\mathbf{1}}$
($(i,u)\in \mathcal{I}_{n+}$).
\end{Definition}

By Lemma \ref{lem:Gx2}, we have the following embedding.
\begin{Theorem}
\label{thm:GTiso}
The ring $\EuScript{T}^{\circ}_{n}(\mathrm{SG})_+$ is isomorphic to
$\mathcal{A}_T(B,x)$ by the correspondence
$T_i(u)\mapsto [\tilde{x}_i(u)]_{\mathbf{1}}$.
\end{Theorem}

The {\em coefficient group $\mathcal{G}(B,y)$
associated with $\mathcal{A}(B,x,y)$}
is the multiplicative subgroup of
the semif\/ield $\mathbb{P}_{\mathrm{univ}}(y)$ generated by all
the coef\/f\/icients $y_{\mathbf{i}}'$ of $\mathcal{A}(B,x,y)$
together with $1+y_{\mathbf{i}}'$.

The following lemma follows from the exchange relation of
coef\/f\/icients \eqref{eq:coef} and the property of the sequence
\eqref{eq:GB2}.

\begin{Lemma}
\label{lem:Gy2}
The family $\{ y_i(u)
\mid (i,u)\in \mathcal{I}'_{n+}\}$
satisfies the Y-system $\mathbb{Y}_{n}(\mathrm{SG})$
by replacing $Y_i(u)$ with $y_i(u)$.
\end{Lemma}

\begin{Example}
Consider the case $n=6$ in Fig.~\ref{fig:labelxG1}.
Let us consider the mutation at the vertex~$(1,1)$ in $Q(0)$,
to which  the variable $y_1(0)$ is attached.
The next time~$(1,1)$ is mutated is in~$Q(10)$,
where $y_1(10)$ is attached.
Meanwhile, between $u=0$ and $u=10$,
the vertices connected to $(1,1)$ in~$Q(u)$ and mutated
are
$(6,2)$ at $u=1,9$, $(6,3)$ at $u=2,8$, $(6,4)$ at $u=3,7$,
$(6,5)$ at $u=4,6$, and  $(6,6)$ and $(6,7)$ at~$u=5$,
Taking account of the directions of the arrows,
we have the relation
\begin{gather*}
y_1(0)y_1(10)=
(1+y_2(1))(1+y_2(9))(1+y_3(2))(1+y_3(8)) (1+y_4(3))\\
\phantom{y_1(0)y_1(10)=}{} \times (1+y_4(7))
(1+ y_5(4))(1+y_5(6)) (1+y_6(5))(1+y_7(5)),
\end{gather*}
which agrees with
\eqref{eq:Y1}.
Similarly, consider the mutation at the vertex $(6,2)$ in $Q(1)$,
to which  the variable $y_2(1)$ is attached.
The next time $(6,2)$ is mutated is in $Q(3)$,
where $y_2(3)$ is attached.
Meanwhile, between $u=1$ and $u=3$,
the vertices connected to $(6,2)$ in $Q(u)$ and mutated
are~$(2,1)$ and~$(6,3)$ at $u=2$.
Taking account of the directions of the arrows,
we have the relation
\[
y_2(1)y_2(3) =
\frac{(1+y_1(2))}{1+y_3(2)^{-1}}.
\]
As the last example, consider the mutation at the vertex $(6,6)$ in $Q(1)$,
to which  the variable~$y_6(1)$ is attached.
The next time $(6,6)$ is mutated is in $Q(3)$,
where $y_6(3)$ is attached.
Meanwhile, between $u=1$ and $u=3$,
the only vertex connected to $(6,6)$ in $Q(u)$ and mutated
is  $(6,5)$ at~$u=2$.
Taking account of the directions of the arrows,
we have the relation
\[
y_6(1)y_6(3) =
\frac{1}{1+y_5(2)^{-1}}.
\]
The other cases can be checked in similar manners.
\end{Example}

\begin{Definition}
The {\em Y-subgroup
$\mathcal{G}_Y(B,y)$
of $\mathcal{G}(B,y)$
associated with the sequence \eqref{eq:Gmutseq}}
is the subgroup of
$\mathcal{G}(B,y)$
generated by
$y_{\mathbf{i}}(u)$
($(\mathbf{i},u)\in \mathbf{I}\times \mathbb{Z}$)
and $1+y_{\mathbf{i}}(u)$
($(\mathbf{i},u):\mathbf{p}_+$ or $\mathbf{p}_-$),
or equivalently,
generated by
$y_i(u)$ and $1+y_i(u)$
($(i,u)\in \mathcal{I}'_{n+}$).
\end{Definition}

By Lemma \ref{lem:Gy2}, we have the following embedding.
\begin{Theorem}
\label{thm:GYiso}
The group $\EuScript{Y}^{\circ}_{n}(\mathrm{SG})_+$ is isomorphic to
$\mathcal{G}_Y(B,y)$ by the correspondence
$Y_i(u)\mapsto y_i(u)$
and $1+Y_i(u)\mapsto 1+y_i(u)$.
\end{Theorem}

\section{Proof of
Theorems \ref{thm:Yperiod}, \ref{thm:dilog},
and \ref{thm:Tperiod}}
\label{sec:proofSG}

In this section we prove
Theorems \ref{thm:Yperiod}, \ref{thm:dilog},
and \ref{thm:Tperiod}
using the method of \cite{Nakanishi09,Inoue10a}.

Let $y=y(0)$ be the initial coef\/f\/icient tuple
of the cluster algebra $\mathcal{A}(B,x,y)$
with $B=B_{n}(\mathrm{SG})$ in the previous section.
Let
$\mathbb{P}_{\mathrm{trop}}(y)$ be
the {\em tropical semifield} for $y$ (Appendix~\ref{sec:groc}(i)).
Let $\pi_{\mathbf{T}}:
\mathbb{P}_{\mathrm{univ}}(y)\rightarrow
\mathbb{P}_{\mathrm{trop}}(y)$, $y_{\mathbf{i}}\mapsto
y_{\mathbf{i}}$ be the projection.
Let $[y_{\mathbf{i}}(u)]_{\mathbf{T}}$
and $[\mathcal{G}_Y(B,y)]_{\mathbf{T}}$
denote the images of $y_{\mathbf{i}}(u)$
and $\mathcal{G}_Y(B,y)$
by the multiplicative group
 homomorphism induced from $\pi_{\mathbf{T}}$, respectively.
They are called the {\em tropical evaluations},
and the resulting relations in
the group $[\mathcal{G}_Y(B,y)]_{\mathbf{T}}$
are called the {\em tropical Y-system}.
They are f\/irst studied in~\cite{Fomin03b}
for cluster algebras of f\/inite types.

We say a (Laurent) monomial $m=\prod_{\mathbf{i}\in \mathbf{I}}
y_{\mathbf{i}}^{k_{\mathbf{i}}}$
is {\em positive} (resp. {\em negative})
if $m\neq 1$ and $k_{\mathbf{i}}\geq 0$
(resp. $k_{\mathbf{i}}\leq 0$)
for any $\mathbf{i}$.
It is known that every monomial $[y_{\mathbf{i}}(u)]_{\mathbf{T}}$
is either positive or negative by \cite[Proposition 5.6]{Fomin07}
and \cite[Theorem 1.7]{Derksen10}.

The next `tropical mutation rule' for $[y_{\mathbf{i}}(u)]_{\mathbf{T}}$
is general and  useful.
\begin{Lemma}
\label{lem:mutation}
Suppose that $y''$ is the coefficient tuple obtained from
the mutation of another coefficient tuple $y'$ at $\mathbf{k}$
with mutation matrix $B'$. Then, for any $\mathbf{i}\neq \mathbf{k}$,
we have the rule:
\begin{enumerate}\itemsep=0pt
\item[$(i)$] $[y''_{\mathbf{i}}]_{\mathbf{T}} =
 [y'_{\mathbf{i}}]_{\mathbf{T}}[y'_{\mathbf{k}}]_{\mathbf{T}}$
 if one of the following conditions holds:
\begin{enumerate}\itemsep=0pt
\item[$(a)$] $B'_{\mathbf{k}\mathbf{i}}> 0$,
 and $[y'_{\mathbf{k}}]_{\mathbf{T}}$ is
positive;
\item[$(b)$]  $B'_{\mathbf{k}\mathbf{i}}< 0$,
 and $[y'_{\mathbf{k}}]_{\mathbf{T}}$ is
negative.
\end{enumerate}
\item[$(ii)$]
$[y''_{\mathbf{i}}]_{\mathbf{T}} =
 [y'_{\mathbf{i}}]_{\mathbf{T}}$
 if one of the following conditions holds:
\begin{enumerate}\itemsep=0pt
\item[$(a)$] $B'_{\mathbf{k}\mathbf{i}}=0$;

\item[$(b)$]
 $B'_{\mathbf{k}\mathbf{i}}> 0$,
 and $[y'_{\mathbf{k}}]_{\mathbf{T}}$ is
 negative;
\item[$(c)$] $B'_{\mathbf{k}\mathbf{i}}< 0$,
 and $[y'_{\mathbf{k}}]_{\mathbf{T}}$ is
 positive.
 \end{enumerate}
 \end{enumerate}
\end{Lemma}
\begin{proof}
This is an immediate consequence of the exchange relation
\eqref{eq:coef} and \eqref{eq:trop}.
\end{proof}

The following properties
of the tropical Y-system are crucial.

\begin{Proposition}
\label{prop:lev2}
 For
$[\mathcal{G}_Y(B,y)]_{\mathbf{T}}$
with $B=B_{n}(\mathrm{SG})$, the following facts hold.

\begin{enumerate}\itemsep=0pt
\item[$(i)$]  For $0 \le u < 4n-2$,
the  monomial $[y_{\mathbf{i}}(u)]_{\mathbf{T}}$
$((\mathbf{i},u):\mathbf{p}_+)$
is negative if and only if $u$ takes the following values.
\[
\begin{cases}
2n-2 \leq u < 4n-2 &
\quad\mbox{\rm for $\mathbf{i}=(1,1),\dots,(n-1,1),(n,2)$},\\
u=2n-2,2n-1, 4n-4, 4n-3 &
\quad\mbox{\rm for $\mathbf{i}=(n,3),\dots,(n,n+1)$}.
\end{cases}
\]
$($Note that for each $\mathbf{i}$, $u$ takes only a part
of the list due to the condition $(\mathbf{i},u):\mathbf{p}_+.)$
\item[$(ii)$]
 We have $[y_{\mathbf{i}}(4n-2)]_{\mathbf{T}}
=y_{\tau^{-1}(\mathbf{i})}$,
where $\tau$ is a bijection $\mathbf{I} \rightarrow
\mathbf{I}$ defined by
\begin{gather*}
(i,1)   \mapsto (\sigma(i),1),
\qquad i=1,\dots,n-1,\\
(n,i')   \mapsto
\begin{cases}
(n,i'),& i'=2,\dots,n-1,\\
(n,n+1),& i'=n,\\
(n,n),& i'=n+1
\end{cases}
\end{gather*}
and $\sigma$ is the permutation in~\eqref{eq:Q2p}.

\item[$(iii)$] The number $N_-$ of the negative monomials
$[y_{\mathbf{i}}(u)]_{\mathbf{T}}$ for $(\mathbf{i},u):\mathbf{p}_+$
in the region $0\leq u < 4n-2$ is $4n-2$.
\end{enumerate}
\end{Proposition}

\begin{proof}
$(i)$ Let us factorize $[y_{\mathbf{i}}(u)]_{\mathbf{T}}
=[y_{\mathbf{i}}(u)]'_{\mathbf{T}}
[y_{\mathbf{i}}(u)]''_{\mathbf{T}}$,
where
 $[y_{\mathbf{i}}(u)]'_{\mathbf{T}}$ is a monomial
in $y_{(i,1)}$ ($i=1,\dots,n-1$)
while
 $[y_{\mathbf{i}}(u)]''_{\mathbf{T}}$ is a monomial
in $y_{(n,i')}$ ($i'=2,\dots,n+1$).
One can independently study
$[y_{\mathbf{i}}(u)]'_{\mathbf{T}}$ and
$[y_{\mathbf{i}}(u)]''_{\mathbf{T}}$.
The claim  (i) follows from the following results,
which  are proved inductively on $u$
by Lemma \ref{lem:mutation} and the results of
\cite{Fomin03b,Fomin07}.

$(a)$ $[y_{\mathbf{i}}(u)]'_{\mathbf{T}}$ part.
This part is easier.
All the monomials $[y_{\mathbf{i}}(u)]'_{\mathbf{T}}$
 which are not 1 for $(\mathbf{i},u):\mathbf{p}_+$
in the region $0 \leq u < 4n-2$ are as follows.
We have $[y_{(i,1)}(2i-2)]'_{\mathbf{T}}=y_{(i,1)}$
$(i=1,\dots,n-1)$,
and also
\begin{alignat}{3}
& [y_{(1,1)}(2n-2)]'_{\mathbf{T}}=y_{(1,1)}^{-1},\qquad &&
[y_{(n,2)}(2n-1)]'_{\mathbf{T}}=y_{(1,1)}^{-1},&\notag\\
& [y_{(2,1)}(2n)]'_{\mathbf{T}}=y_{(1,1)}^{-1}y_{(2,1)}^{-1},\qquad &&
[y_{(n,2)}(2n+1)]'_{\mathbf{T}}=y_{(2,1)}^{-1},&\notag\\
\label{eq:y'part}
& [y_{(3,1)}(2n+2)]'_{\mathbf{T}}=y_{(2,1)}^{-1}y_{(3,1)}^{-1},\qquad &&
[y_{(n,2)}(2n+3)]'_{\mathbf{T}}=y_{(3,1)}^{-1},&\\
&\qquad \vdots & &\qquad \vdots & \notag\\
& [y_{(n-1,1)}(4n-6)]'_{\mathbf{T}}=y_{(n-2,1)}^{-1}y_{(n-1,1)}^{-1},\qquad &&
[y_{(n,2)}(4n-5)]'_{\mathbf{T}}=y_{(n-1,1)}^{-1},&\notag\\
& [y_{(1,1)}(4n-4)]'_{\mathbf{T}}=y_{(n-1,1)}^{-1}.&&& \notag
\end{alignat}

$(b)$  $[y_{\mathbf{i}}(u)]''_{\mathbf{T}}$ part.
Below we list all the monomials $[y_{\mathbf{i}}(u)]''_{\mathbf{T}}$
which are not 1 for $(\mathbf{i},u):\mathbf{p}_+$
in the region $0 \leq u < 4n-2$.
We separate the region $0 \leq u < 4n-2$ into four parts
corresponding to the decomposition
$4n-2=(2n-2)+2+(2n-4)+2$,
where $2n-2$ and $2n-4$ are the Coxeter numbers
of $D_{n}$ and $D_{n-1}$.

Region I: $0\leq u < 2n-2$.
All the  monomials
$[y_{(n,i')}(u)]''_{\mathbf{T}}$ ($i'=2,\dots,n+1$)
 for $(\mathbf{(n,i')},u):\mathbf{p}_+$
are identif\/ied with
the {\em positive roots of $D_n$} as in
\cite[Proposition 10.7]{Fomin07};
therefore, they are positive.
Here, $D_{n}$ is identif\/ied with the subgraph of $X_n$
consisting of vertices $2, \dots, n+1$.

Region II: $u=2n-2$, $2n-1$.
We have, for even $n$,
\begin{gather*}
[y_{(n,i')}(2n-2)]''_{\mathbf{T}} =y_{(n,i')}^{-1},\qquad  i'=3,5,\dots,n-1 ,\nonumber\\
[y_{(n,i')}(2n-1)]''_{\mathbf{T}} =y_{(n,i')}^{-1},\qquad  i'=2,4,\dots,n,n+1 ,%\label{eq:yt1}
\end{gather*}
and, for odd $n$,
\begin{gather*}
[y_{(n,i')}(2n-2)]''_{\mathbf{T}} =y_{(n,i')}^{-1}, \qquad  i'=3,5,\dots,n-2 ,\nonumber\\
[y_{(n,n)}(2n-2)]''_{\mathbf{T}} =y_{(n,n+1)}^{-1},
\qquad [y_{(n,n+1)}(2n-2)]''_{\mathbf{T}}=y_{(n,n)}^{-1},\nonumber\\
[y_{(n,i')}(2n-1)]''_{\mathbf{T}} =y_{(n,i')}^{-1},\qquad i'=2,4,\dots,n-1.%\label{eq:yt4}
\end{gather*}

Region III: $2n \leq u < 4n-4$.
All the  monomials
$[y_{(n,i')} (u)]''_{\mathbf{T}}$ ($i'=3,\dots,n+1$)
 for $((n,i'),u):\mathbf{p}_+$
are  identif\/ied with
the {\em positive roots of $D_{n-1}$},
therefore, they are positive.
Here, $D_{n-1}$ is identif\/ied with the subgraph of $X_n$
consisting of vertices $3, \dots, n+1$.

Region IV: $u=4n-4$, $4n-3$. We have, for even $n$,
\begin{gather}
[y_{(n,i')}(4n-4)]''_{\mathbf{T}} =y_{(n,i')}^{-1}, \qquad i'=3,5,\dots,n-1,\nonumber\\
[y_{(n,i')}(4n-3)]''_{\mathbf{T}} =y_{(n,i')}^{-1}, \qquad i'=2,4,\dots,n-2,\nonumber\\
[y_{(n,n)}(4n-3)]''_{\mathbf{T}} =y_{(n,n+1)}^{-1},
\qquad [y_{(n,n+1)}(4n-3)]''_{\mathbf{T}}=y_{(n,n)}^{-1},\label{eq:yt3}
\end{gather}
and, for odd $n$,
\begin{gather}
[y_{(n,i')}(4n-4)]''_{\mathbf{T}} =y_{(n,i')}^{-1}, \qquad  i'=3,5,\dots,n-2 ,\nonumber\\
[y_{(n,n)}(4n-4)]''_{\mathbf{T}} =y_{(n,n+1)}^{-1},
\qquad [y_{(n,n+1)}(4n-4)]''_{\mathbf{T}}=y_{(n,n)}^{-1},\nonumber\\
[y_{(n,i')}(4n-3)]''_{\mathbf{T}} =y_{(n,i')}^{-1},\qquad i'=2,4,\dots,n-1.\label{eq:yt6}
\end{gather}

Besides, we have the sequences of  monomials
which appear over Regions~III and~IV;
for even~$n$,
\begin{gather}
 [y_{(2,1)}(2n)]''_{\mathbf{T}}=y_{(n,2)}^{-1}y_{(n,3)}^{-1},
\qquad [y_{(3,1)}(2n+2)]''_{\mathbf{T}}=y_{(n,4)}^{-1}y_{(n,5)}^{-1},
\qquad \dots,\nonumber\\
 [y_{(n/2,1)}(3n-4)]''_{\mathbf{T}}=y_{(n,n-2)}^{-1}y_{(n,n-1)}^{-1},\qquad
 [y_{(n/2+1,1)}(3n-2)]''_{\mathbf{T}}=y_{(n,n-1)}^{-1}y_{(n,n)}^{-1}y_{(n,n+1)}^{-1},\nonumber\\
 [y_{(n/2+2,1)}(3n)]''_{\mathbf{T}}=y_{(n,n-3)}^{-1}y_{(n,n-2)}^{-1},\qquad
\dots,\nonumber\\
  [y_{(n-1,1)}(4n-6)]''_{\mathbf{T}}=y_{(n,3)}^{-1}y_{(n,4)}^{-1},
\qquad [y_{(1,1)}(4n-4)]''_{\mathbf{T}}=y_{(n,2)}^{-1},\label{eq:yt2}
\end{gather}
and, for odd $n$, the middle three terms are replaced with
\begin{gather}
 [y_{((n-1)/2,1)}(3n-5)]''_{\mathbf{T}}=y_{(n,n-3)}^{-1}y_{(n,n-2)}^{-1},\nonumber\\
 [y_{((n+1)/2,1)}(3n-3)]''_{\mathbf{T}}=
y_{(n,n-1)}^{-1}y_{(n,n)}^{-1}y_{(n,n+1)}^{-1},\nonumber\\
 [y_{(n+3)/2,1)}(3n-1)]''_{\mathbf{T}}=y_{(n,n-2)}^{-1}y_{(n,n-1)}^{-1}.\label{eq:yt5}
\end{gather}

$(ii)$ They  follow from
\eqref{eq:yt3}--\eqref{eq:yt5}.

$(iii)$ By $(i)$, for each  $\mathbf{i}$
the numbers of the negative monomials
$[y_{\mathbf{i}}(u)]_{\mathrm{T}}$
in the region is~2 for $\mathbf{i}= (1,1)$,
 1 for $\mathbf{i}= (i,1)$ ($i=2,\dots,n-1$),
$n$ for $\mathbf{i}= (n,2)$.
and $2$ for $\mathbf{i}=(n,i')$ ($i'=3,\dots,n+1$).
Summing up, we have $N_-=4n-2$.
\end{proof}

Now we prove
Theorems \ref{thm:Yperiod} and \ref{thm:Tperiod}.
It follows from  a very general theorem
\cite[Theorem 5.1]{Inoue10a}
(based on the work by Plamondon \cite{Plamondon10a,Plamondon10b})
that the cluster variables
$x_{\mathbf{i}}(u)$ and coef\/f\/icients $y_{\mathbf{i}}(u)$
have the same periodicity with
$[y_{\mathbf{i}}(u)]_{\mathbf{T}}$ as in Proposition \ref{prop:lev2}$(ii)$,
namely,
$x_{\mathbf{i}}(4n-2)=x_{\tau^{-1}(\mathbf{i})}$
and
$y_{\mathbf{i}}(4n-2)=y_{\tau^{-1}(\mathbf{i})}$.
It follows that, under the labelling introduced in Lemma
\ref{lem:gmap}, we have
\begin{gather}
\tilde{x}_i(u+4n-2) =\tilde{x}_{\omega(i)}(u),
\qquad (i,u):\mathbf{P}_+,\nonumber\\
y_i(u+4n-2) =y_{\omega(i)}(u),
\qquad (i,u):\mathbf{P}'_+,\label{eq:xperiod}
\end{gather}
where $\omega$ is the one in Theorem~\ref{thm:Yperiod}$(i)$.
Then, thanks to the isomorphisms
in Theorems~\ref{thm:GTiso} and~\ref{thm:GYiso},
and also by the isomorphisms
$\EuScript{T}^{\circ}_{n}(\mathrm{SG})_+
\simeq
\EuScript{T}^{\circ}_{n}(\mathrm{SG})_-
$
and
$\EuScript{Y}^{\circ}_{n}(\mathrm{SG})_+
\simeq
\EuScript{Y}^{\circ}_{n}(\mathrm{SG})_-
$,
we obtain
Theorems~\ref{thm:Yperiod} and~\ref{thm:Tperiod}.

Since the {\em $F$-polynomials\/}
(Appendix \ref{sec:groc}(vii))
are def\/ined as a certain specialization of cluster variables,
they satisfy the same periodicity as the cluster variables.
Let $\tilde{F}_i(u)$ ($(i,u):\mathbf{P}_+$) be the $F$-polynomial
for $\tilde{x}_i(u)$.
Then, from \eqref{eq:xperiod} we have
\begin{gather}
\label{eq:F}
\tilde{F}_i(u+4n-2)=\tilde{F}_{\omega(i)}(u),
\qquad (i,u):\mathbf{P}_+.
\end{gather}

Next we prove
Theorem \ref{thm:dilog}.
Let $\bigwedge^2 \mathbb{P}_{\mathrm{univ}}(y)$
be the quotient of
the additive Abelian  group
$\mathbb{P}_{\mathrm{univ}}(y)\otimes_{\mathbb{Z}}
\mathbb{P}_{\mathrm{univ}}(y)$
by the subgroup generated by
symmetric tensors~\cite{Frenkel95,Chapoton05}.
\begin{Lemma}
\label{lem:dilog}
In $\bigwedge^2 \mathbb{P}_{\mathrm{univ}}(y)$, we have
\[
\sum_{
\genfrac{}{}{0pt}{1}
{
(i,u):\mathbf{P}'_+
}
{
0\leq u < 4n-2
}
} y_i(u)\wedge (1+y_i(u))=0.
\]
\end{Lemma}
\begin{proof}
One can prove it in a similar way as in \cite{Nakanishi09,Inoue10a},
using \cite[Proposition 3.13]{Fomin07} and~\eqref{eq:F}.$\!\!\!\!\!$
\end{proof}

Applying the method of \cite{Frenkel95,Chapoton05,Nakanishi09},
we immediately obtain the following theorem from
Lemma~\ref{lem:dilog} and Proposition~\ref{prop:lev2}$(iii)$.

\begin{Theorem}
For any semifield homomorphism $\varphi:\mathbb{P}_{\mathrm{univ}}(y)
\rightarrow \mathbb{R}_+$, we have the following identity.
\begin{gather*}
%\label{eq:DI2}
\frac{6}{\pi^2}
\sum_{
\genfrac{}{}{0pt}{1}
{
(i,u):\mathbf{P}'_+
}
{
0\leq u < 4n-2
}
}
L\left(
\frac{\varphi(y_i(u))}{1+\varphi(y_i(u))}
\right)
 =4n-2.
\end{gather*}
\end{Theorem}
This is equivalent to Theorem \ref{thm:dilog}.

\section{Cluster algebras for RSG Y-systems}
\label{sec:clusterRSG}

In this section we identify $\mathbb{T}_{n}(\mathrm{RSG})$ and
$\mathbb{Y}_{n}(\mathrm{RSG})$ as relations for cluster variables
and coef\/f\/i\-cients of the cluster algebra
associated with a certain quiver $Q_{n}(\mathrm{RSG})$.
For those  things which are quite parallel to the SG case,
we skip their precise descriptions just by saying `as before'
unless they are not obvious.

\subsection{Parity decompositions of T and Y-systems}

For $(i,u)\in \tilde{\mathcal{I}}_{n}$,
we set the parity condition $\mathbf{P}_{+}$ by,
for even $n$,
\begin{gather*}
%\label{eq:GPcond1r}
\mathbf{P}_{+}: \
\mbox{$i+u$ is even}, \quad i=1,\dots,n-2,
\end{gather*}
and, for odd $n$,
\begin{gather*}
%\label{eq:GPcond2r}
\mathbf{P}_{+}: \
\begin{cases}
\mbox{$u$ is even},& i=1,\\
\mbox{$i+u$ is even},& i=2,\dots,n-2.
\end{cases}
\end{gather*}
Def\/ine
$\tilde{\mathcal{I}}_{n\varepsilon}$
and $\EuScript{T}^{\circ}_{n}(\mathrm{RSG})_{\varepsilon}$
($\varepsilon=\pm$)
as before.
Then, we have
$\EuScript{T}^{\circ}_{n}(\mathrm{RSG})_+
\simeq
\EuScript{T}^{\circ}_{n}(\mathrm{RSG})_-
$
by $T_i(u)\mapsto T_i(u+1)$ and
\begin{gather*}
%\label{eq:Tfactr}
\EuScript{T}^{\circ}_{n}(\mathrm{RSG})
\simeq
\EuScript{T}^{\circ}_{n}(\mathrm{RSG})_+
\otimes_{\mathbb{Z}}
\EuScript{T}^{\circ}_{n}(\mathrm{RSG})_-.
\end{gather*}

For  $(i,u)\in \tilde{\mathcal{I}}_{n}$,
we set the parity condition $\mathbf{P}'_{+}$ by
\begin{gather*}
%\label{eq:GP'cond1r}
\mathbf{P}'_{+}: \
\mbox{$i+u$ is odd},\quad i=1,\dots,n-2.
\end{gather*}
We have
\[
(i,u): \ \mathbf{P}'_+ \ \Longleftrightarrow\
\textstyle (i,u\pm d_i):\ \mathbf{P}_+,
\]
where $d_i$ is given in \eqref{eq:di}.
Def\/ine $\tilde{\mathcal{I}}'_{n\varepsilon}$
and $\EuScript{Y}^{\circ}_{n}(\mathrm{RSG})_{\varepsilon}$
($\varepsilon=\pm$) as before.
Then, we have
$\EuScript{Y}^{\circ}_{n}(\mathrm{RSG})_+
\simeq
\EuScript{Y}^{\circ}_{n}(\mathrm{RSG})_-
$
by $Y_i(u)\mapsto Y_i(u+1)$,
$1+Y_i(u)\mapsto 1+Y_i(u+1)$,
 and
\begin{gather*}
%\label{eq:Yfactr}
\EuScript{Y}^{\circ}_{n}(\mathrm{RSG})
\simeq
\EuScript{Y}^{\circ}_{n}(\mathrm{RSG})_+
\times
\EuScript{Y}^{\circ}_{n}(\mathrm{RSG})_-.
\end{gather*}

From now on, we mainly treat the $+$ parts,
$\EuScript{T}^{\circ}_{n}(\mathrm{RSG})_+$
and
$\EuScript{Y}^{\circ}_{n}(\mathrm{RSG})_+$.

\subsection[Quiver $Q_n(\mathrm{RSG})$]{Quiver $\boldsymbol{Q_n(\mathrm{RSG})}$}
\label{subsect:quiverr}

%%%%%%%%%%%%%%%%%%%%%%%%%%%%%
% SG
\begin{figure}[t]
\centerline{\includegraphics{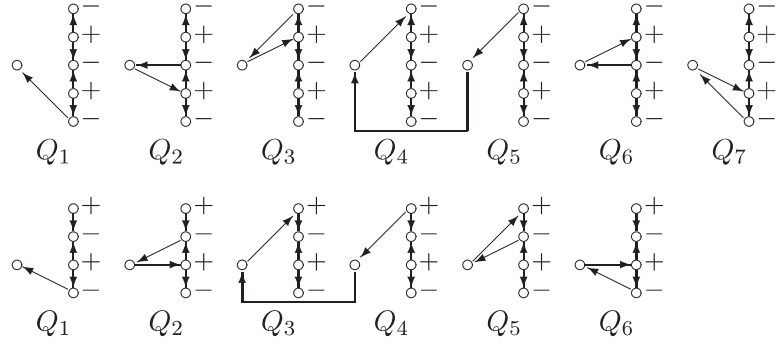}}

\caption{The quiver $Q_{n}(\mathrm{RSG})$  for $n=8$
(upper) and for $n=7$ (lower),
where, except for the leftmost vertex of each quiver $Q_i$,
 all the  vertices in the same position
 in  $n-1$ quivers $Q_1,\dots,Q_{n-1}$ are
identif\/ied.
 We have an arrow between the leftmost vertices
of $Q_4$ and $Q_5$ for $n=8$, and $Q_3$ and $Q_4$ for $n=7$.}
\label{fig:quiverGr}
\end{figure}

With each $n \geq 4$ we associate
 a quiver $Q_{n}(\mathrm{RSG})$ as below.
First,  the cases $n=8$ and $n=7$ are given in
 Fig.~\ref{fig:quiverGr},
where, except for the leftmost vertex of each quiver $Q_i$,
 all the  vertices in the same position
 in  $n-1$ quivers $Q_1$,\dots,$Q_{n-1}$ are
identif\/ied.
For a general $n$, the quiver
$Q_{n}(\mathrm{RSG})$ is def\/ined by naturally extending
these examples;
in particular, we have an arrow  from the leftmost vertex
of $Q_{n/2+1}$ to the leftmost vertex of $Q_{n/2}$
for even $n$,
and
an arrow  from the leftmost vertex
of $Q_{(n+1)/2}$ to the leftmost vertex of $Q_{(n-1)/2}$
for odd $n$.
Also we assign the property $+/-$
to  each vertex, except for the leftmost one in each $Q_i$,
 as in Fig.~\ref{fig:quiverGr}.

Let us choose  the index set $\tilde{\mathbf{I}}$
of the vertices of $Q_{n}(\mathrm{RSG})$
naturally obtained from
the
index set~$\mathbf{I}$
of the vertices of $Q_{n}(\mathrm{SG})$
by the restriction to the
vertices of~$Q_{n}(\mathrm{RSG})$.
Thus, $i=1,\dots,n$; and $i'=1$ if $i\neq n$
and $i'=2,\dots,n-2$ if $i=n$.
Then, we def\/ine $\tilde{\mathbf{I}}_+$, $\tilde{\mathbf{I}}_-$,
$\mu_+$, $\mu_-$,
and $\tilde{w}(Q_{n}(\mathrm{RSG}))$
as before.

\begin{Lemma}
\label{lem:GQmutr}
Let $Q(0):=Q_{n}(\mathrm{RSG})$.
We have the periodic sequence of mutations of qui\-vers~\eqref{eq:GB2},
where the quiver $Q(2p)$ $(p=1,\dots,n-2)$ is given by
\eqref{eq:Q2p}.
\end{Lemma}

\begin{Example}
The mutation sequence
 \eqref{eq:GB2} for $n=6$ is explicitly given
in Fig.~\ref{fig:labelxG1r}.
Let us explain why we need the `extra arrow' from
$Q_4$ to $Q_3$ in $Q(0)$ in this example.
Suppose that we do not have it.
Then, we will have an extra arrow from
$Q_3$ to $Q_4$ in $Q(2)$ due to the mutation at $u=1$.
This extra arrow will still remain in $Q(4)$,
where the vertex $(3,1)$ will be mutated.
This conf\/licts with $\mathbb{T}_n(\mathrm{RSG})$
and $\mathbb{Y}_n(\mathrm{RSG})$ eventually.
The extra arrow is, thus, necessary
 as a~precaution to avoid this conf\/liction.
\begin{figure}[t]
\centerline{\includegraphics{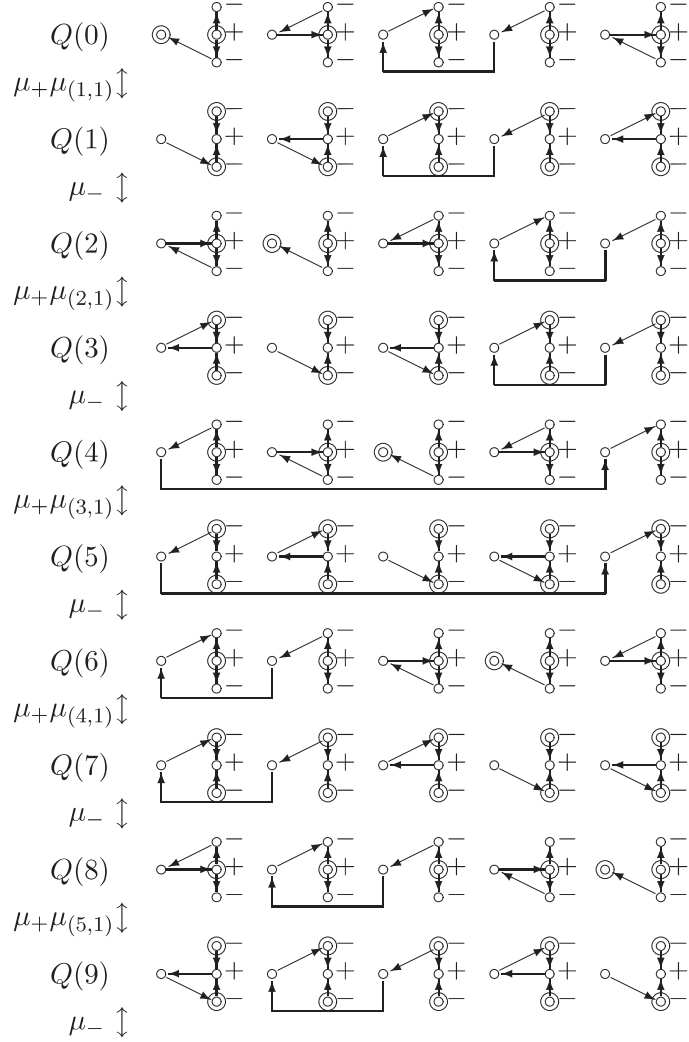}}

\caption{The mutation sequence of the quiver $Q_{n}(\mathrm{RSG})$
in  \eqref{eq:GB2} for $n=6$.
The encircled vertices correspond to
the mutation points $(\mathbf{i},u): \mathbf{p}_+$
in the forward direction.}
\label{fig:labelxG1r}
\end{figure}
\end{Example}

\subsection{Embedding maps}
Let
 $B=B_{n}(\mathrm{RSG})$
be the skew-symmetric matrix corresponding
to the quiver $Q_{n}(\mathrm{RSG})$.
Let $\mathcal{A}(B,x,y)$
be the cluster algebra
 with coef\/f\/icients
in the universal
semif\/ield
$\mathbb{P}_{\mathrm{univ}}(y)$.

In view of Lemma \ref{lem:GQmutr}
we set $x(0)=x$, $y(0)=y$ and def\/ine
clusters $x(u)=(x_{\mathbf{i}}(u))_{\mathbf{i}\in \tilde{\mathbf{I}}}$
 ($u\in \mathbb{Z}$)
 and coef\/f\/icient tuples $y(u)=(y_\mathbf{i}(u))_{\mathbf{i}\in
 \tilde{\mathbf{I}}}$
 ($u\in \mathbb{Z}$)
by the sequence of mutations \eqref{eq:Gmutseq} as before.

{\samepage For  $(\mathbf{i},u)\in
 \tilde{\mathbf{I}}\times \mathbb{Z}$,
we set the parity condition $\mathbf{p}_+$ by the same rule
as \eqref{eq:GQparity}, where
$\mathbf{I}_{\pm}$ is replaced with $\tilde{\mathbf{I}}_{\pm}$.
We def\/ine the condition $\mathbf{p}_-$
by $(\mathbf{i},u):\mathbf{p}_- \Longleftrightarrow
(\mathbf{i},u-1):\mathbf{p}_+$.

}

\begin{Lemma}
\label{lem:gmapr}
Below $\equiv$ means the equivalence modulo $(2n-2)\mathbb{Z}$.

\begin{enumerate}\itemsep=0pt
\item[$(i)$] The map
\begin{gather*}
g:  \ \tilde{\mathcal{I}}_{n+} \rightarrow
\{ (\mathbf{i},u)\in \tilde{\mathbf{I}}\times \mathbb{Z} \mid
 (\mathbf{i},u): \mathbf{p}_+
\},\\
\phantom{g:  \ }{} \ (i,u-d_i) \mapsto
\begin{cases}
((j+1,1),u), & i= 1; \
u\equiv 2j
\ (j=0,1,\dots,n-2),\\
((n,i),u), & i=2,\dots,n-2
\end{cases}
\end{gather*}
is a bijection.

\item[$(ii)$]
The map
\begin{gather*}
g': \ \tilde{\mathcal{I}}'_{n+} \rightarrow
\{ (\mathbf{i},u)\in \tilde{\mathbf{I}}\times \mathbb{Z} \mid
 (\mathbf{i},u): \mathbf{p}_+
\},\\
\phantom{g': \ }{} \ (i,u) \mapsto
\begin{cases}
((j+1,1),u), & i= 1; \
u\equiv 2j
\ (j=0,1,\dots,n-2), \\
((n,i),u), & i=2,\dots,n-2
\end{cases}
\end{gather*}
is a bijection.
\end{enumerate}
\end{Lemma}

Based on Lemma~\ref{lem:gmapr},
we introduce alternative notations
$\tilde{x}_i(u-d_i):=x_{\mathbf{i}}(u)$
for $(i,u-d_i)\in \tilde{\mathcal{I}}_{n+}$
with  $(\mathbf{i},u)=g((i,u-d_i))$
and
$y_i(u):=y_{\mathbf{i}}(u)$
for $(i,u)\in \tilde{\mathcal{I}}'_{n+}$
with $(\mathbf{i},u)=g'((i,u))$,
respectively.

Let $\mathcal{A}(B,x)$ be the cluster algebra
with trivial coef\/f\/icients.

\begin{Definition}
\label{def:T-subr}
The {\em T-subalgebra
$\mathcal{A}_T(B,x)$
of ${\mathcal{A}}(B,x)$
associated with the sequence \eqref{eq:Gmutseq}}
is the subring of
${\mathcal{A}}(B,x)$
generated by
$[x_{\mathbf{i}}(u)]_{\mathbf{1}}$
($(\mathbf{i},u)\in \tilde{\mathbf{I}}\times \mathbb{Z}$),
or equivalently,
generated by
$[\tilde{x}_i(u)]_{\mathbf{1}}$
($(i,u)\in \tilde{\mathcal{I}}_{n+}$).
\end{Definition}

By the lemma parallel
to Lemma \ref{lem:Gx2}, we have the following embedding.
\begin{Theorem}
\label{thm:GTisor}
The ring $\EuScript{T}^{\circ}_{n}(\mathrm{RSG})_+$ is isomorphic to
$\mathcal{A}_T(B,x)$ by the correspondence
$T_i(u)\mapsto [\tilde{x}_i(u)]_{\mathbf{1}}$.
\end{Theorem}

\begin{Example}
Consider the case $n=6$ in Fig.~\ref{fig:labelxG1r}.
Consider the mutation at the vertex $(6,4)$ in $Q(1)$,
to which  the variable $\tilde{x}_4(0)$ is attached.
The next time $(6,4)$ is mutated is in $Q(3)$,
where $\tilde{x}_4(2)$ is attached.
Meanwhile, the vertices connected to $(6,4)$ in $Q(1)$
are $(3,1)$, $(4,1)$, $(5,1)$ and $(6,3)$,
and the  variable attached to them are equal to
$\tilde{x}_1(-1)$, $\tilde{x}_1(1)$, $\tilde{x}_1(3)$, and $\tilde{x}_3(1)$, respectively.
Taking account of the directions of the arrows,
we have the relation
\[
\tilde{x}_4(0)\tilde{x}_4(2)=
\frac{y_4(1)}{1+y_4(1)} \tilde{x}_1(-1)\tilde{x}_1(3)
+ \frac{1}{1+y_4(1)} \tilde{x}_1(1)\tilde{x}_3(1),
\]
which agree with~\eqref{eq:T1r}.
The other cases are the same as in the SG case.
\end{Example}

Let  $\mathcal{G}(B,y)$ be the
coef\/f\/icient group
associated with $\mathcal{A}(B,x,y)$ as before.

\begin{Definition}
The {\em Y-subgroup
$\mathcal{G}_Y(B,y)$
of $\mathcal{G}(B,y)$
associated with the sequence~\eqref{eq:Gmutseq}}
is the subgroup of
$\mathcal{G}(B,y)$
generated by
$y_{\mathbf{i}}(u)$
($(\mathbf{i},u)\in \tilde{\mathbf{I}}\times \mathbb{Z}$)
and $1+y_{\mathbf{i}}(u)$
($(\mathbf{i},u):\mathbf{p}_+$ or $\mathbf{p}_-$),
or equivalently,
generated by $y_i(u)$ and $1+y_i(u)$
$((i,u)\in \tilde{\mathcal{I}}'_{n+})$.
\end{Definition}

By the lemma parallel to Lemma \ref{lem:Gy2},
 we have the following embedding.
\begin{Theorem}
\label{thm:GYisor}
The group $\EuScript{Y}^{\circ}_{n}(\mathrm{RSG})_+$ is isomorphic to
$\mathcal{G}_Y(B,y)$ by the correspondence
$Y_i(u)\mapsto y_i(u)$
and $1+Y_i(u)\mapsto 1+y_i(u)$.
\end{Theorem}

\begin{Example}
Consider the case $n=6$ in Fig.~\ref{fig:labelxG1r}.
Let us consider the mutation at the vertex $(1,1)$ in $Q(0)$,
to which  the variable $y_1(0)$ is attached.
The next time $(1,1)$ is mutated is in $Q(10)$,
where $y_1(10)$ is attached.
Meanwhile, between $u=0$ and $u=10$,
the vertices connected to $(1,1)$ in $Q(u)$ and mutated
are
$(6,2)$ at $u=1,9$, $(6,3)$ at $u=2,8$, and $(6,4)$ at $u=3,5,7$.
Taking account of the directions of the arrows,
we have the relation
\[
y_1(0)y_1(10)=
\frac{
(1+y_2(1))(1+y_2(9))(1+y_3(2))(1+y_3(8)) (1+y_4(3))
 (1+y_4(7))
}
{1+y_4(5)^{-1}},
\]
which agrees with
\eqref{eq:Y1r}.
The other cases are the same as in the SG case.
\end{Example}

\section{Proof of
Theorems \ref{thm:Yperiodr}, \ref{thm:dilogr},
and \ref{thm:Tperiodr}}
\label{sec:proofRSG}

In this section we prove
Theorems \ref{thm:Yperiodr}, \ref{thm:dilogr},
and \ref{thm:Tperiodr}
by the same method as in Section \ref{sec:proofSG}.
\label{sec:proofSGr}

Let $B_n(\mathrm{RSG})$ be the one in the previous section.
The following properties
of the tropical Y-system are crucial.

\begin{Proposition}
\label{prop:lev2r}
 For
$[\mathcal{G}_Y(B,y)]_{\mathbf{T}}$
with $B=B_{n}(\mathrm{RSG})$, the following facts hold.

\begin{enumerate}\itemsep=0pt
\item[$(i)$]  For $0 \le u < 4n-2$,
the  monomial $[y_{\mathbf{i}}(u)]_{\mathbf{T}}$
$((\mathbf{i},u):\mathbf{p}_+)$
is negative if and only if $u$ takes the following values.
\begin{gather*}
\begin{cases}
2n-2 \leq u < 4n-2  &
\quad\mbox{\rm for $\mathbf{i}=(1,1),\dots,(n-1,1)$},\\
u=n-2,n-1; 2n-2\leq u < 4n-2 &
\quad\mbox{\rm for $\mathbf{i}=(n,2)$},\\
u=n-2,n-1,2n-2,2n-1, &
\\
\quad\quad3n-3, 3n-2,4n-4,4n-3 &\quad\mbox{\rm for $\mathbf{i}=(n,3),\dots,(n,n-2)$}.
\end{cases}
\end{gather*}
$($Note that for each $\mathbf{i}$, $u$ takes only a part
of the list due to the condition $(\mathbf{i},u):\mathbf{p}_+.)$

\item[$(ii)$]
 We have $[y_{\mathbf{i}}(4n-2)]_{\mathbf{T}}
=y_{\tau^{-1}(\mathbf{i})}$,
where $\tau$ is a bijection $\mathbf{I} \rightarrow
\mathbf{I}$ defined by
\begin{gather*}
(i,1)  \mapsto (\sigma(i),1),
\qquad i=1,\dots,n-1,\\
(n,i')   \mapsto
(n,i'),\qquad i'=2,\dots,n-2,
\end{gather*}
and $\sigma$ is the permutation in \eqref{eq:Q2p}.

\item[$(iii)$] The number $N_-$ of the negative monomials
$[y_{\mathbf{i}}(u)]_{\mathbf{T}}$ for $(\mathbf{i},u):\mathbf{p}_+$
in the region $0\leq u < 4n-2$ is $6n-15$.
\end{enumerate}
\end{Proposition}

\begin{proof}
$(i)$ Let us factorize $[y_{\mathbf{i}}(u)]_{\mathbf{T}}
=[y_{\mathbf{i}}(u)]'_{\mathbf{T}}
[y_{\mathbf{i}}(u)]''_{\mathbf{T}}$,
where
 $[y_{\mathbf{i}}(u)]'_{\mathbf{T}}$ is a monomial
in $y_{(i,1)}$ ($i=1,\dots,n-1$)
while
 $[y_{\mathbf{i}}(u)]''_{\mathbf{T}}$ is a monomial
in $y_{(n,i')}$ ($i'=2,\dots,n-2$).
The claim $(i)$ follows from the following results.

$(a)$ $[y_{\mathbf{i}}(u)]'_{\mathbf{T}}$ part.
This part is exactly the same as
\eqref{eq:y'part} in the $B_n(\mathrm{SG})$ case.

$(b)$  $[y_{\mathbf{i}}(u)]''_{\mathbf{T}}$ part.
Below we list all the monomials $[y_{\mathbf{i}}(u)]''_{\mathbf{T}}$
which are not 1 for $(\mathbf{i},u):\mathbf{p}_+$
in the region $0 \leq u < 4n-2$.
We separate the region $0 \leq u < 4n-2$ into eight parts
corresponding to the decomposition
$4n-2=(n-2)+2+(n-2)+2+(n-3)+2+(n-3)+2$,
where $n-2$ and $n-3$ are the Coxeter numbers
of $A_{n-3}$ and $A_{n-4}$.

Region I: $0\leq u < n-2$.
All the  monomials
$[y_{(n,i')}(u)]''_{\mathbf{T}}$ ($i'=2,\dots,n-2$)
 for $((n,i'),u):\mathbf{p}_+$
are identif\/ied with
the positive roots of $A_{n-3}$;
therefore, they are positive.
Here, $A_{n-3}$ is identif\/ied with  the subgraph of $X_n$
consisting of vertices $2, \dots, n-2$.

Region II: $u=n-2$, $n-1$.
We have, for even $n$,
\begin{gather*}
[y_{(n,i')}(n-2)]''_{\mathbf{T}} =y_{(n,n-i')}^{-1}, \qquad i'=3,5,\dots,n-3,\\
[y_{(n,i')}(n-1)]''_{\mathbf{T}} =y_{(n,n-i')}^{-1},\qquad  i'=2,4,\dots,n-2,
\end{gather*}
and, for odd $n$,
\begin{gather*}
[y_{(n,i')}(n-2)]''_{\mathbf{T}} =y_{(n,n-i')}^{-1}, \qquad i'=2,4,\dots,n-3,\\
[y_{(n,i')}(n-1)]''_{\mathbf{T}} =y_{(n,n-i')}^{-1},\qquad i'=3,5,\dots,n-2.
\end{gather*}

Region III: $n \leq u < 2n-2$.
Once again, all the  monomials
$[y_{(n,i')} (u)]''_{\mathbf{T}}$ ($i'=2,\dots,n-2$)
 for $((n,i'),u):\mathbf{p}_+$
are  identif\/ied with
the positive roots of $A_{n-3}$;
therefore, they are positive.

Region IV: $n=2n-2, 2n-1$. We have
\begin{gather*}
%\label{eq:yt1r}
[y_{(n,i')}(2n-2)]''_{\mathbf{T}} =y_{(n,i')}^{-1}, \qquad i'=3,5,\dots,\\
[y_{(n,i')}(2n-1)]''_{\mathbf{T}} =y_{(n,i')}^{-1}, \qquad  i'=2,4,\dots.
\end{gather*}

Region V: $2n \leq u < 3n-3$.
All the  monomials
$[y_{(n,i')}(u)]''_{\mathbf{T}}$ ($i'=3,\dots,n-2$)
 for $((n,i'),u):\mathbf{p}_+$
are identif\/ied with
the positive roots of $A_{n-4}$;
therefore, they are positive.
Here, $A_{n-4}$ is identif\/ied with the subgraph of $X_n$
consisting of vertices $3, \dots, n-2$.

Region VI: $u=3n-3, 3n-2$.
We have, for even $n$,
\begin{gather*}
[y_{(n,i')}(3n-3)]''_{\mathbf{T}} =y_{(n,n+1-i')}^{-1},\qquad  i'=4,6,\dots,n-2,\nonumber\\
[y_{(n,i')}(3n-2)]''_{\mathbf{T}} =y_{(n,n+1-i')}^{-1},\qquad  i'=3,5,\dots,n-3,%\label{eq:yt1r}
\end{gather*}
and, for odd $n$,
\begin{gather*}
[y_{(n,i')}(3n-3)]''_{\mathbf{T}} =y_{(n,n+1-i')}^{-1},\qquad i'=3,5,\dots,n-2,\nonumber\\
[y_{(n,i')}(3n-2)]''_{\mathbf{T}} =y_{(n,n+1-i')}^{-1},\qquad  i'=4,6,\dots,n-3.%\label{eq:yt4r}
\end{gather*}

Region VII: $3n-1 \leq u < 4n-4$.
Once again, all the  monomials
$[y_{(n,i')}(u)]''_{\mathbf{T}}$ ($i'=3,\dots,n-2$)
 for $((n,i'),u):\mathbf{p}_+$
are identif\/ied with
the positive roots of $A_{n-4}$;
therefore, they are positive.

Region VIII: $u=4n-4, 4n-3$. We have
\begin{gather}
[y_{(n,i')}(4n-4)]''_{\mathbf{T}} =y_{(n,i')}^{-1}, \qquad i'=3,5,\dots,\nonumber\\
[y_{(n,i')}(4n-3)]''_{\mathbf{T}} =y_{(n,i')}^{-1}, \qquad i'=2,4,\dots.\label{eq:yt3r}
\end{gather}

Besides, we have the sequences of  monomials
which appear over Regions V--VIII;
for even $n$,
\begin{gather}
 [y_{(2,1)}(2n)]''_{\mathbf{T}}=y_{(n,2)}^{-1}y_{(n,3)}^{-1},
\qquad y_{(3,1)}(2n+2)]''_{\mathbf{T}}=y_{(n,4)}^{-1}y_{(n,5)}^{-1},\qquad
\dots,\nonumber\\
 [y_{(n/2,1)}(3n-4)]''_{\mathbf{T}}=y_{(n,n-2)}^{-1},\qquad
[y_{(n/2+2,1)}(3n)]''_{\mathbf{T}}=y_{(n,n-3)}^{-1}y_{(n,n-2)}^{-1},\qquad
\dots,\nonumber\\
 [y_{(n-1,1)}(4n-6)]''_{\mathbf{T}}=y_{(n,3)}^{-1}y_{(n,4)}^{-1},
\qquad [y_{(1,1)}(4n-4)]''_{\mathbf{T}}=y_{(n,2)}^{-1},\label{eq:yt2r}
\end{gather}
and, for odd $n$, the middle two terms are replaced with
\begin{gather}
[y_{((n-1)/2,1)}(3n-5)]''_{\mathbf{T}}=y_{(n,n-3)}^{-1}y_{(n,n-2)}^{-1},\nonumber\\
[y_{(n+3)/2,1)}(3n-1)]''_{\mathbf{T}}=y_{(n,n-2)}^{-1}.\label{eq:yt5r}
\end{gather}

$(ii)$ They  follow from \eqref{eq:yt3r}--\eqref{eq:yt5r}.

$(iii)$ By $(i)$, for each  $\mathbf{i}$
the numbers of the negative monomials
$[y_{\mathbf{i}}(u)]_{\mathrm{T}}$
in the region is
 2 for $\mathbf{i}= (1,1)$,
 1 for $\mathbf{i}= (i,1)$ ($i=2,\dots,n-1$),
$n+1$ for $\mathbf{i}= (n,2)$.
and $4$ for $\mathbf{i}= (n,i')$ ($i'=3,\dots,n-2$).
Summing up, we have $N_-=6n-15$.
\end{proof}

Theorems \ref{thm:Yperiodr}, \ref{thm:Tperiodr}
and \ref{thm:dilogr} follow from Proposition
\ref{prop:lev2r} as before.

\section{Further extension}
\label{sec:extension}

Seeing that the cluster algebraic
setting perfectly works for the case~\eqref{eq:case1},
it is rather natural to expect that
the method in this paper is,
at least in principle, applicable for a general rational~$\xi$.
On the other hand, working out in full generality seems
a complicated task, and we do not pursue it in this paper.
However, repeating our method to a little more general
case
\begin{gather}
\label{eq:case2}
\xi=\frac{n-1}{mn-m+1}=
\cfrac{1}{ m +
\cfrac{1}{n-1}
}
\end{gather}
is not dif\/f\/icult.
Since the proof is mostly the repetition of the previous case $m=1$,
we concentrate on exhibiting the results.

\subsection{Further result for SG Y-systems}

\begin{figure}[t]
\centerline{\includegraphics{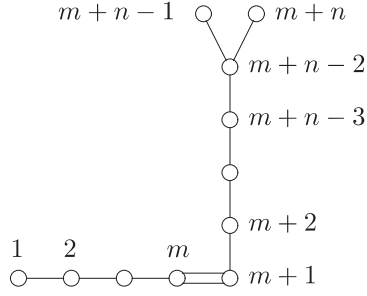}}

\caption{The diagram $X_{m,n}$.}
\label{fig:Xmn}
\end{figure}

With a pair of  integers $m \geq 1$ and  $n \geq 4$, we associate a
diagram $X_{m,n}$ in Fig.~\ref{fig:Xmn}, following~\cite{Tateo95}.
Let $\mathcal{I}_{m,n}=\{1,\dots,m+n\}\times \mathbb{Z}$.

\begin{Definition}[\cite{Tateo95}]
Fix a pair of integers $m \geq 1$ and $n \geq 4$.
The sine-Gordon (SG)
Y-system $\mathbb{Y}_{m,n}(\mathrm{SG})$
is the following system of relations for
a family of variables
 $\{Y_i(u) \mid (i,u)\in \mathcal{I}_{m,n} \}$,
\begin{gather*}
Y_i(u-n+1)Y_i(u+n-1)
 =\prod_{j:j\sim i} (1+Y_j(u)),
\qquad i=1,\dots,m-1,
\nonumber\\
Y_m\left(u-n+1\right)
Y_m\left(u+n-1\right)
 = (1+Y_{m-1}(u))
\nonumber\\
\qquad {}\times
\Biggl(\prod_{j=m+1}^{m+n-2} (1+Y_j(u- m-n+1+j))
(1+Y_j(u+m+n-1-j))\Biggr)\nonumber\\
\qquad {} \times (1+Y_{m+n-1} (u))(1+Y_{m+n}(u)),
\nonumber\\
Y_{m+1}(u-1) Y_{m+1}(u+1)
=
\frac{1+Y_m(u)}{1+Y_{m+2}(u)^{-1}}
,\nonumber\\
Y_i(u-1)Y_i(u+1)
=
\frac{1}{\prod\limits_{j:j\sim i} (1+Y_j(u)^{-1})},
\qquad i=m+2,\dots,m+n,%\label{eq:Y1m}
\end{gather*}
where $j\sim i$ means that $j$ is adjacent to $i$ in $X_{m,n}$.
\end{Definition}

\begin{Definition}
Fix a pair of integers $m \geq 1$ and $n \geq 4$.
The sine-Gordon (SG)
T-system $\mathbb{T}_{m,n}(\mathrm{SG})$
is the following system of relations for
a family of variables
 $\{T_i(u) \mid (i,u)\in \mathcal{I}_{m,n} \}$,
\begin{gather*}
T_i(u-n+1)T_i(u+n-1)
=
\prod_{j:j\sim i} T_j(u)+1,
\qquad i=1,\dots,m,\nonumber\\
T_{m+1}(u-1) T_{m+1}(u+1)
=
T_m(u-n+2)T_m(u+n-2)+T_{m+2}(u),
\nonumber\\
T_i(u-1)T_i(u+1)=T_m(u-m-n+1+i)\nonumber
\\
\phantom{T_i(u-1)T_i(u+1)=}{}\times T_m(u+m+n-1-i)
+ \!\prod_{j: j\sim i} T_j(u),\! \qquad i=m+2,\dots,m+n-2,
\nonumber\\
T_{m+n-1}(u-1)T_{m+n-1}(u+1)
=
T_m(u) + T_{m+n-2}(u),
\nonumber\\
T_{m+n}(u-1)T_{m+n}(u+1)
=
T_m(u) + T_{m+n-2}(u),%\label{eq:T1m}
\end{gather*}
where $j\sim i$ means that $j$ is adjacent to $i$ in $X_{m,n}$.
\end{Definition}

We def\/ine the associated multiplicative Abelian
group $\EuScript{Y}^{\circ}_{m,n}(\mathrm{SG})$
and the commutative ring
$\EuScript{T}^{\circ}_{m,n}(\mathrm{SG})$
 as before.

\begin{Theorem}\label{thm:Yperiodm}\qquad
\begin{enumerate}\itemsep=0pt
\item[$(i)$] {\rm (Conjectured by \cite{Tateo95}.)}
The following relations hold in $\EuScript{Y}^{\circ}_{m,n}(\mathrm{SG})$.
\begin{itemize}\itemsep=0pt
\item[$(a)$] Suppose that both $m$ and $n$ are even.
$($In other words, $mn-m+n$ is even.$)$
Then, we have the periodicity: $Y_i(u+2(mn-m+n))=Y_i(u)$.
$($There is no half periodicity.$)$
\item[$(b)$] Suppose that at least one of $m$ or $n$ is odd.
$($In other words, $mn-m+n$ is odd.$)$
Then, we have the half periodicity: $Y_i(u+2(mn-m+n))=Y_{\omega(i)}(u)$, where
$\omega$ is an involution of the set $\{1,\dots,m+n\}$
defined by $\omega(m+n-1)=m+n$, $\omega(m+n)=m+n-1$,
and $\omega(i)=i$ $(i=1,\dots,m+n-2)$.
Therefore, we have the full periodicity: $Y_i(u+4(mn-m+n))=Y_i(u)$.
\end{itemize}

\item[$(ii)$]  The same periodicity holds
in  $\EuScript{T}^{\circ}_{m,n}(\mathrm{SG})$
by replacing $Y_i(u)$ in $(i)$ with $T_i(u)$.

\item[$(iii)$] {\rm (Conjectured by \cite{Tateo95}.)}
Suppose that a family of positive real numbers
$\{ Y_i(u) \mid (i,u)\in \mathcal{I}_n\}$ satisfies
$\mathbb{Y}_{m,n}(\mathrm{SG})$.
Then, we have the identities
\begin{gather*}
%\label{eq:DIm}
\frac{6}{\pi^2}
\sum_{
\genfrac{}{}{0pt}{1}
{
(i,u)\in \mathcal{I}_{m,n}
}
{
0\leq u < 4(mn-m+n)
}
}
L\left(
\frac{Y_i(u)}{1+Y_i(u)}
\right)
=4(m+1)(mn-m+n),\\
%\label{eq:DI'm}
\frac{6}{\pi^2}
\sum_{
\genfrac{}{}{0pt}{1}
{
(i,u)\in \mathcal{I}_{m,n}
}
{
0\leq u < 4(mn-m+n)
}
}
L\left(
\frac{1}{1+Y_i(u)}
\right)
=4(n-1)(mn-m+n).
\end{gather*}
\end{enumerate}
\end{Theorem}

In our proof of Theorem \ref{thm:Yperiodm}
we have  a natural interpretation
of the  full/half period
\[
2(mn-m+n)=h(D_n)+2+ m(h(D_{n-1})+2).
\]

\begin{Remark}
Actually,
 $\mathbb{Y}_{m,n}(\mathrm{SG})$ and  $\mathbb{T}_{m,n}(\mathrm{SG})$
are also considered for $n=3$,
and they coincide with the Y and T-systems of
type $B_{m+1}$ with level $2$ in \cite{Kuniba94a}.
Theorem~\ref{thm:Yperiodm} remains valid for $n=3$ due to~\cite{Inoue10a}.
\end{Remark}

\begin{figure}[t]
\centerline{\includegraphics{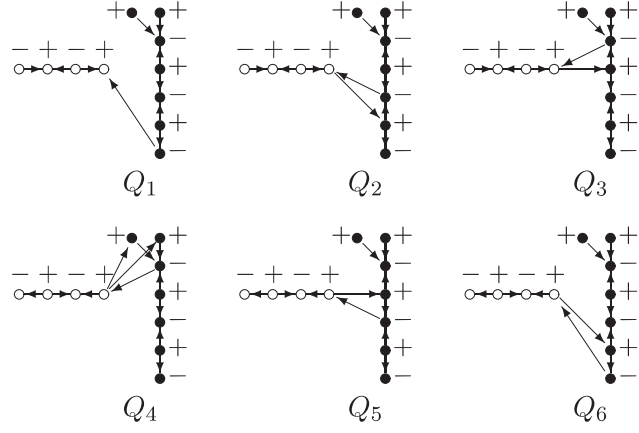}}

\caption{The quiver $Q_{m,n}(\mathrm{SG})$
for $(m,n)=(4,7)$,
where,
 all the  vertices with $\bullet$ in the same position
 in  the quivers $Q_1,\dots,Q_{6}$ are
identif\/ied.
Note that the arrows between the vertices with $\circ$
in $Q_1$, $Q_2$, $Q_3$ and $Q_4$, $Q_5$, $Q_6$ are opposite.}
\label{fig:quiverGm}
\end{figure}

The cluster algebraic formulation
of $\mathbb{Y}_{m,n}(\mathrm{SG})$
and $\mathbb{T}_{m,n}(\mathrm{SG})$ is done
by the quiver
$Q_{m,n}(\mathrm{SG})$ def\/ined as follows.
As a rather general example, the case $(m,n)=(4,7)$
is given in
 Fig.~\ref{fig:quiverGm},
where,  all the  vertices with $\bullet$ in the same position
 in  the quivers $Q_1,\dots,Q_{6}$ are
identif\/ied.
Note that the arrows between the vertices with $\circ$
in $Q_1$, $Q_2$, $Q_3$ and $Q_4$, $Q_5$, $Q_6$ are opposite.
For a general $n$,  the quiver $Q_{n}(\mathrm{SG})$ is def\/ined
by naturally extending this example.
Namely, add $m-1$ vertices to the left of the leftmost
vertex   in each quiver $Q_i$ of
$Q_{n}(\mathrm{SG})$.
We assign the property $+/-$
to  each vertex as in Fig.~\ref{fig:quiverGm}
such that the leftmost vertex of each $Q_i$ has the property $-$ for even $m$ and
$+$ for odd $m$.
Then, we put the arrows between the vertices with $\circ$
such that, for even $n$, each vertex with property $(\circ,-)$
is a source for $Q_1,\dots, Q_{n/2}$ and a sink for
$Q_{n/2+1},\dots, Q_{n-1}$;
and, for odd $n$, each vertex with property $(\circ,-)$
is a source for $Q_1,\dots, Q_{(n-1)/2}$ and a sink for
$Q_{(n+1)/2+1},\dots, Q_{n-1}$.

Let us choose  the index set $\mathbf{I}$
of the vertices of $Q_{m,n}(\mathrm{SG})$
so that,
under the  natural
identif\/ication with $X_{m,n}$,
 $\mathbf{i}=(i,i')\in \mathbf{I}$
($i=1,\dots,n-1$; $i'=1,\dots,m$) represents
the vertex  with $\circ$ in~$Q_i$,
and $\mathbf{i}=(n,i')\in \mathbf{I}$
($i'=m+1,\dots,m+n$)
 represents
the vertex  with $\bullet$ in any quiver $Q_i$.

Let $\mathbf{I}_+^{\bullet}$
(resp.\ $\mathbf{I}_-^{\bullet}$)
denote the set of the vertices with
property  $(\bullet,+)$ (resp.\  $(\bullet,-)$).
Let $\mathbf{I}_{+,i}^{\circ}$
(resp.\ $\mathbf{I}_{-,i}^{\circ}$)
denote the set of the vertices with
property  $(\circ,+)$ (resp.\  $(\circ,-)$) in the $i$th quiver $Q_i$.

We def\/ine composite mutations,
\begin{gather}
\label{eq:mupm2m}
\mu_+^{\bullet}=\prod_{\mathbf{i}\in\mathbf{I}_+^{\bullet}}
\mu_{\mathbf{i}},
\qquad
\mu_-^{\bullet}=\prod_{\mathbf{i}\in\mathbf{I}_-^{\bullet}}
\mu_{\mathbf{i}},
\qquad
\mu_{+,i}^{\circ}=\prod_{\mathbf{i}\in\mathbf{I}_{+,i}^{\circ}}
\mu_{\mathbf{i}},
\qquad
\mu_{-,i}^{\circ}=\prod_{\mathbf{i}\in\mathbf{I}_{-,i}^{\circ}}
\mu_{\mathbf{i}}.
\end{gather}

\begin{Lemma}
\label{lem:GQmutm}
Let $Q(0):=Q_{m,n}(\mathrm{SG})$.
We have the following periodic sequence of mutations of quivers:
for even $n$,
\begin{gather}
Q(0)
\mathop{\longleftrightarrow}^{{\mu_+^{\bullet}} \mu_{+,1}^{\circ}}
 Q(1)
\mathop{\longleftrightarrow}^{\mu_-^{\bullet}\mu_{-,n/2+1}^{\circ}}
Q(2)
\mathop{\longleftrightarrow}^{\mu_+^{\bullet} \mu_{+,2}^{\circ}}
Q(3)
\mathop{\longleftrightarrow}^{\mu_-^{\bullet}\mu_{-,n/2+2}^{\circ}}
\ \cdots\
\nonumber\\
\phantom{Q(0)}{}
\mathop{\longleftrightarrow}^{\mu_+^{\bullet}\mu_{+,n/2}^{\circ}}
Q(n-1)
\mathop{\longleftrightarrow}^{{\mu_-^{\bullet}} \mu_{-,1}^{\circ}}
 Q(n)
\mathop{\longleftrightarrow}^{\mu_+^{\bullet}\mu_{+,n/2+1}^{\circ}}
 Q(n+1)
\mathop{\longleftrightarrow}^{\mu_-^{\bullet}\mu_{-,2}^{\circ}}
\ \cdots\
\nonumber\\
\phantom{Q(0)}{}
\mathop{\longleftrightarrow}^{\mu_-^{\bullet}\mu_{-,n/2-1}^{\circ}}
Q(2n-4)
\mathop{\longleftrightarrow}^{\mu_+^{\bullet}\mu_{+,n-1}^{\circ}}
Q(2n-3)
\mathop{\longleftrightarrow}^{\mu_-^{\bullet}\mu_{-,n/2}^{\circ}}
Q(2n-2)=Q(0),\label{eq:GB2m}
\end{gather}
and, for odd $n$,
\begin{gather}
Q(0)
\mathop{\longleftrightarrow}^{{\mu_+^{\bullet}} \mu_{+,1}^{\circ}
\mu_{-,(n+1)/2}^{\circ}}
 Q(1)
\mathop{\longleftrightarrow}^{\mu_-^{\bullet}}
Q(2)
\mathop{\longleftrightarrow}^{\mu_+^{\bullet} \mu_{+,2}^{\circ}
\mu_{-,(n+3)/2}^{\circ}}
Q(3)
\mathop{\longleftrightarrow}^{\mu_-^{\bullet}}
\ \cdots\
\nonumber\\
\phantom{Q(0)}{}
\mathop{\longleftrightarrow}^{\mu_-^{\bullet}}
Q(n-1)
\mathop{\longleftrightarrow}^{{\mu_+^{\bullet}}
\mu_{+,(n+1)/2}^{\circ} \mu_{-,1}^{\circ}
}
 Q(n)
\mathop{\longleftrightarrow}^{\mu_-^{\bullet}}
Q(n+1)
\mathop{\longleftrightarrow}^{{\mu_+^{\bullet}}
\mu_{+,(n+3)/2}^{\circ} \mu_{-,2}^{\circ}
}
\ \cdots
\nonumber\\
\phantom{Q(0)}{}
\mathop{\longleftrightarrow}^{\mu_-^{\bullet}}
Q(2n-4)
\mathop{\longleftrightarrow}^{\mu_+^{\bullet}\mu_{+,n-1}^{\circ}
\mu_{-,(n-1)/2}^{\circ}}
Q(2n-3)
\mathop{\longleftrightarrow}^{\mu_-^{\bullet}}
Q(2n-2)=Q(0).\label{eq:GB2mm}
\end{gather}
\end{Lemma}

For $(\mathbf{i},u)\in \mathbf{I}\times \mathbb{Z}$,
write $(\mathbf{i},u):\mathbf{p}_+$
if $(\mathbf{i},u)$ is the forward mutation points
in \eqref{eq:GB2m} for even~$n$ and  in \eqref{eq:GB2mm} for odd $n$,
modulo $(2n-2)\mathbb{Z}$ for $u$.
Then, one can repeat and extend all the arguments for $m=1$
in Sections~\ref{sec:clusterSG} and~\ref{sec:proofSG},
prove the following proposition, and obtain
 Theorem~\ref{thm:Yperiodm}.

\begin{Proposition}
\label{prop:lev2m}
 For
$[\mathcal{G}_Y(B,y)]_{\mathbf{T}}$
with $B=B_{m,n}(\mathrm{SG})$, the following facts hold.
\begin{enumerate}\itemsep=0pt
\item[$(i)$]  For $0 \le u < 2(mn-m+n)$,
the  monomial $[y_{\mathbf{i}}(u)]_{\mathbf{T}}$
$((\mathbf{i},u):\mathbf{p}_+)$
is negative if and only if $u$ takes the values
$2n-2 \leq u < 2(mn-m+n)$
for $\mathbf{i}=(i,i')$ $(i=1,\dots,n-1;i'=1,\dots,m)$ and
for $\mathbf{i}
=(n,m+1)$,
and
$u=2k(n-1),2k(n-1)+1\ (k=1,\dots,m+1)$
for $\mathbf{i}=(n,m+2),\dots,(n,m+n)$.
\item[$(ii)$]
 We have $[y_{\mathbf{i}}(2(mn-m+n))]_{\mathbf{T}}
=y_{\tau^{-1}(\mathbf{i})}$,
where $\tau$ is a bijection $\mathbf{I} \rightarrow
\mathbf{I}$ defined as follows.
If both $m$ and $n$ are even,
\begin{gather}
(i,i')   \mapsto (\sigma(i),i'),
\qquad i=1,\dots,n-1,\quad i'=1,\dots,m,\nonumber\\
(n,i')   \mapsto
(n,i'),\qquad i'=m+1,\dots,m+n\label{eq:tau1}
\end{gather}
and $\sigma$ is the permutation in~\eqref{eq:Q2p}.
If at least one of
$m$ or $n$ is odd,
we modify $\tau$ in~\eqref{eq:tau1} by
$(n,m+n-1)\mapsto (n,m+n)$ and $(n,m+n)\mapsto (n,m+n-1)$.
$($The rest are the same as in~\eqref{eq:tau1}.$)$

\item[$(iii)$] The number $N_-$ of the negative monomials
$[y_{\mathbf{i}}(u)]_{\mathbf{T}}$ for $(\mathbf{i},u):\mathbf{p}_+$
in the region $0\leq u < 2(mn-m+n)$ is $(m+1)(mn-m+n)$.
\end{enumerate}
\end{Proposition}

\subsection{Further result for RSG Y-system}

The RSG case is quite similar.
Let $\tilde{\mathcal{I}}_{m,n}=\{1,\dots,m+n-3\}\times \mathbb{Z}$.

\begin{Definition}[\cite{Tateo95}]
Fix a pair of integers $m \geq 1$ and $n \geq 4$.
The reduced sine-Gordon (RSG)
Y-system $\mathbb{Y}_{m,n}(\mathrm{RSG})$
is the following system of relations for
a family of variables
 $\{Y_i(u) \mid (i,u)\in \tilde{\mathcal{I}}_{m,n} \}$,
\begin{gather*}
Y_i(u-n+1)Y_i(u+n-1)
=\prod_{j:j\sim i} (1+Y_j(u)),
\qquad i=1,\dots,m-1,
\nonumber\\
Y_m\left(u-n+1\right)
Y_m\left(u+n-1\right)
= (1+Y_{m-1}(u))
\nonumber\\
\qquad {}\times
\Biggl(\prod_{j=m+1}^{m+n-3} (1+Y_j(u-m-n+1+j))
(1+Y_j(u+m+n-1-j))\Biggr)\nonumber\\
\qquad{} \times (1+Y_{m+n-3} (u)^{-1})^{-1},
\nonumber\\
Y_{m+1}(u-1) Y_{m+1}(u+1)
=
\frac
{1+Y_m(u)}
{1+Y_{m+2}(u)^{-1}}
,\nonumber\\
Y_i(u-1)Y_i(u+1)
=
\frac{1}
{\prod\limits_{j:j\sim i} (1+Y_j(u)^{-1})
},
\qquad i=m+2,\dots,m+n-3,%\label{eq:Y1mr}
\end{gather*}
where $j\sim i$ means that $j\leq m+n-3$ is adjacent to $i$ in $X_{m,n}$.
\end{Definition}

\begin{Definition}
Fix a pair of integers $m \geq 1$ and $n \geq 4$.
The reduced sine-Gordon (RSG)
T-system $\mathbb{T}_{m,n}(\mathrm{RSG})$
is the following system of relations for
a family of variables
 $\{T_i(u) \mid (i,u)\in \tilde{\mathcal{I}}_{m,n} \}$,
\begin{gather*}
T_i(u-n+1)T_i(u+n-1)
=
\prod_{j:j\sim i} T_j(u)+1,
\qquad i=1,\dots,m,\nonumber\\
T_{m+1}(u-1) T_{m+1}(u+1)
=
T_m(u-n+2)T_m(u+n-2)+T_{m+2}(u),
\nonumber\\
T_i(u-1)T_i(u+1)= T_m(u-m-n+1+i)\nonumber\\
\phantom{T_i(u-1)T_i(u+1)=}{}\times T_m(u+m+n-1-i)
+ \!\prod_{j: j\sim i} T_j(u), \!\qquad i=m+2,\dots,m+n-4,
\nonumber\\
T_{m+n-3}(u-1)T_{m+n-3}(u+1)
=
T_{m}(u-2)T_{m}(u+2) +T_{m}(u) T_{m+n-4}(u),%\label{eq:T1mr}
\end{gather*}
where $j\sim i$ means that $j\leq m+n-3$ is adjacent to $i$ in $X_{m,n}$.
\end{Definition}

We def\/ine the associated multiplicative Abelian
group $\EuScript{Y}^{\circ}_{m,n}(\mathrm{RSG})$
and the commutative ring
$\EuScript{T}^{\circ}_{m,n}(\mathrm{RSG})$
 as before.

\begin{Theorem}\qquad\label{thm:Yperiodmr}
\begin{enumerate}\itemsep=0pt
\item[$(i)$] {\rm (Conjectured by \cite{Tateo95}.)}
The following relations hold in $\EuScript{Y}^{\circ}_{m,n}(\mathrm{RSG})$.\\
Periodicity: $Y_i(u+2(mn-m+n))=Y_i(u)$.

\item[$(ii)$] %(Conjectured by \cite{Tateo95}.)
The following relations hold in $\EuScript{T}^{\circ}_{m,n}(\mathrm{RSG})$.
\\
Periodicity: $T_i(u+2(mn-m+n))=T_i(u)$.

\item[$(iii)$] {\rm  (Conjectured by \cite{Tateo95}.)}
Suppose that a family of positive real numbers
$\{ Y_i(u) \mid (i,u)\in \mathcal{I}_n\}$ satisfies
$\mathbb{Y}_{m,n}(\mathrm{RSG})$.
Then, we have the identities
\begin{gather*}
%\label{eq:DImr}
\frac{6}{\pi^2}
\sum_{
\genfrac{}{}{0pt}{1}
{
(i,u)\in \mathcal{I}_{m,n}
}
{
0\leq u < 2(mn-m+n)
}
}
L\left(
\frac{Y_i(u)}{1+Y_i(u)}
\right)
 =2\big(nm^2-m^2+3mn-8m+ 2n-6\big),\\
%\label{eq:DI'mr}
\frac{6}{\pi^2}
\sum_{
\genfrac{}{}{0pt}{1}
{
(i,u)\in \mathcal{I}_{m,n}
}
{
0\leq u < 2(mn-m+n)
}
}
L\left(
\frac{1}{1+Y_i(u)}
\right)
 =2\big(n^2m-6nm+11m+n^2-5n+6\big).
\end{gather*}
\end{enumerate}
\end{Theorem}

In our proof of Theorem~\ref{thm:Yperiodm}
we have  a natural interpretation
of the  period
\[
2(mn-m+n)=2\{h(A_{n-3})+2+ m(h(A_{n-4})+2)\}.
\]

\begin{figure}[t]
\centerline{\includegraphics{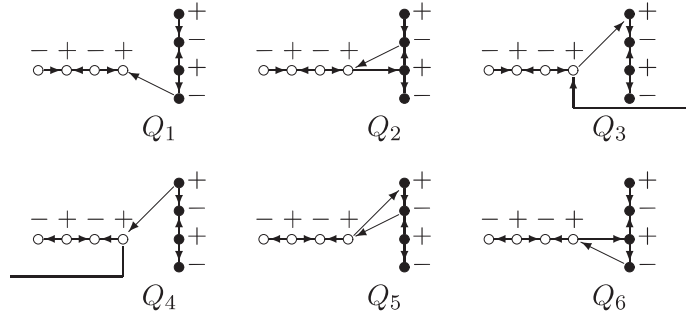}}
\caption{The quiver $Q_{m,n}(\mathrm{RSG})$
for $(m,n)=(4,7)$,
where,
 all the  vertices with $\bullet$ in the same position
 in  the quivers $Q_1,\dots,Q_{6}$ are
identif\/ied.
 We have an arrow between the fourth vertices (from the left)
of~$Q_3$ and~$Q_4$.}
\label{fig:quiverGmr}
\end{figure}

The cluster algebraic formulation
of $\mathbb{Y}_{m,n}(\mathrm{RSG})$
and $\mathbb{T}_{m,n}(\mathrm{RSG})$ is done
by the quiver
$Q_{m,n}(\mathrm{RSG})$ def\/ined as follows.
Add $m-1$ vertices to the left of the leftmost
vertex   in each quiver $Q_i$ of
$Q_{n}(\mathrm{RSG})$
and put the arrows between them  exactly in
the same way as in the quiver $Q_{m,n}(\mathrm{SG})$.
As an example, the case $(m,n)=(4,7)$
is given in
 Fig.~\ref{fig:quiverGmr}.
We also assign the properties $+/-$ and $\circ/\bullet$
to  each vertex
in a similar way to $Q_{m,n}(\mathrm{SG})$
 as in Fig.~\ref{fig:quiverGmr}.

Let us choose  the index set $\tilde{\mathbf{I}}$
of the vertices of $Q_{m,n}(\mathrm{RSG})$
naturally obtained from
the
index set $\mathbf{I}$
of the vertices of $Q_{m,n}(\mathrm{SG})$
by the restriction to the
vertices of $Q_{m,n}(\mathrm{RSG})$.
Thus, $\mathbf{i}=(i,i')\in \tilde{\mathbf{I}}$
($i=1,\dots,n-1$; $i'=1,\dots,m$) represents
the vertex  with $\circ$ in $Q_i$,
and $\mathbf{i}=(n,i')\in \tilde{\mathbf{I}}$
($i'=m+1,\dots,m+n-3$)
 represents
the vertex  with $\bullet$ in any quiver $Q_i$.
Then, we def\/ine
 $\tilde{\mathbf{I}}_+^{\bullet}$,
 $\tilde{\mathbf{I}}_-^{\bullet}$,
 $\tilde{\mathbf{I}}_{+,i}^{\circ}$,
 $\tilde{\mathbf{I}}_{-,i}^{\circ}$
and the corresponding mutations by \eqref{eq:mupm2m}.

\begin{Lemma}
\label{lem:GQmutmr}
Let $Q(0):=Q_{m,n}(\mathrm{RSG})$.
We have the periodic sequence of mutations of qui\-vers~\eqref{eq:GB2m} for even $n$ and
quivers \eqref{eq:GB2mm} for odd $n$.
\end{Lemma}

For $(\mathbf{i},u)\in \mathbf{I}\times \mathbb{Z}$,
write $(\mathbf{i},u):\mathbf{p}_+$
if $(\mathbf{i},u)$ is the forward mutation points
in \eqref{eq:GB2m} for even $n$ and  in \eqref{eq:GB2mm} for odd $n$,
modulo $(2n-2)\mathbb{Z}$ for $u$.
Then, one can repeat and extend all the arguments for $m=1$
in Sections~\ref{sec:clusterRSG} and~\ref{sec:proofRSG},
prove the following proposition, and obtain Theorem~\ref{thm:Yperiodmr}.

\begin{Proposition}
\label{prop:lev2mr}
 For
$[\mathcal{G}_Y(B,y)]_{\mathbf{T}}$
with $B=B_{m,n}(\mathrm{RSG})$, the following facts hold.
\begin{enumerate}\itemsep=0pt
\item[$(i)$]  For $0 \le u < 2(mn-m+n)$,
the  monomial $[y_{\mathbf{i}}(u)]_{\mathbf{T}}$
$((\mathbf{i},u):\mathbf{p}_+)$
is negative if and only if $u$ takes the values
$2n-2 \leq u < 2(mn-m+n)$
for $\mathbf{i}=(i,i')$ $(i=1,\dots,n-1;i'=1,\dots,m)$,
$u=n-2,n-1$, $2n-2 \leq u < 2(mn-m+n)$
for $\mathbf{i}=(n,m+1)$,
and
$u=n-2,n-1,k(n-1),k(n-1)+1\ (k=2,\dots,2m+2)$
for $\mathbf{i}=(n,m+2),\dots,(n,m+n-3)$.

\item[$(ii)$]
 We have $[y_{\mathbf{i}}(2(mn-m+n))]_{\mathbf{T}}
=y_{\tau^{-1}(\mathbf{i})}$,
where $\tau$ is a bijection $\tilde{\mathbf{I}} \rightarrow
\tilde{\mathbf{I}}$ defined by
\begin{gather*}
(i,i')   \mapsto (\sigma(i),i'),
\qquad i=1,\dots,n-1,\quad   i'=1,\dots,m,\\
(n,i')  \mapsto
(n,i'), \qquad i'=m+1,\dots,m+n-3
\end{gather*}
and $\sigma$ is the permutation in \eqref{eq:Q2p}.

\item[$(iii)$] The number $N_-$ of the negative monomials
$[y_{\mathbf{i}}(u)]_{\mathbf{T}}$ for $(\mathbf{i},u):\mathbf{p}_+$
in the region $0\leq u < 2(mn-m+n)$ is $nm^2-m^2+3mn-8m+ 2n-6$.
\end{enumerate}
\end{Proposition}

\subsection{Remarks by referee}
\label{subsec:mutseq}

The content of this subsection is based
on the communication from the referee of this paper.

Recall that two quivers are said to be {\em mutation-equivalent\/}
if one is obtained from the other by a sequence of  mutations.

The following important observation was made by the referee.

{\bf Fact 1.}
 The quiver $Q_{m,n}(\mathrm{SG})$ in
 Fig.~\ref{fig:quiverGm} (Fig.~\ref{fig:quiverG} for $m=1$)
is mutation-equivalent to the quiver of type $D_{mn-m+n}$,
i.e., its underlying graph is the Dynkin diagram of type
$D_{mn-m+n}$.
Similarly,
the quiver $Q_{m,n}(\mathrm{RSG})$ in
Fig.~\ref{fig:quiverGmr} (Fig.~\ref{fig:quiverGr}
for $m=1$)
is mutation-equivalent to the quiver of type $A_{mn-m+n-3}$.
{\em In other words,
the cluster algebra $\mathcal{A}(B,x,y)$
with $B=B_{m,n}(\mathrm{SG})$ $($resp.\ $B=B_{m,n}(\mathrm{RSG}))$
 is the cluster algebra
of type $D_{mn-m+n}$ $($resp. type $A_{mn-m+n-3})$.}

\begin{figure}[t]
\centerline{\includegraphics{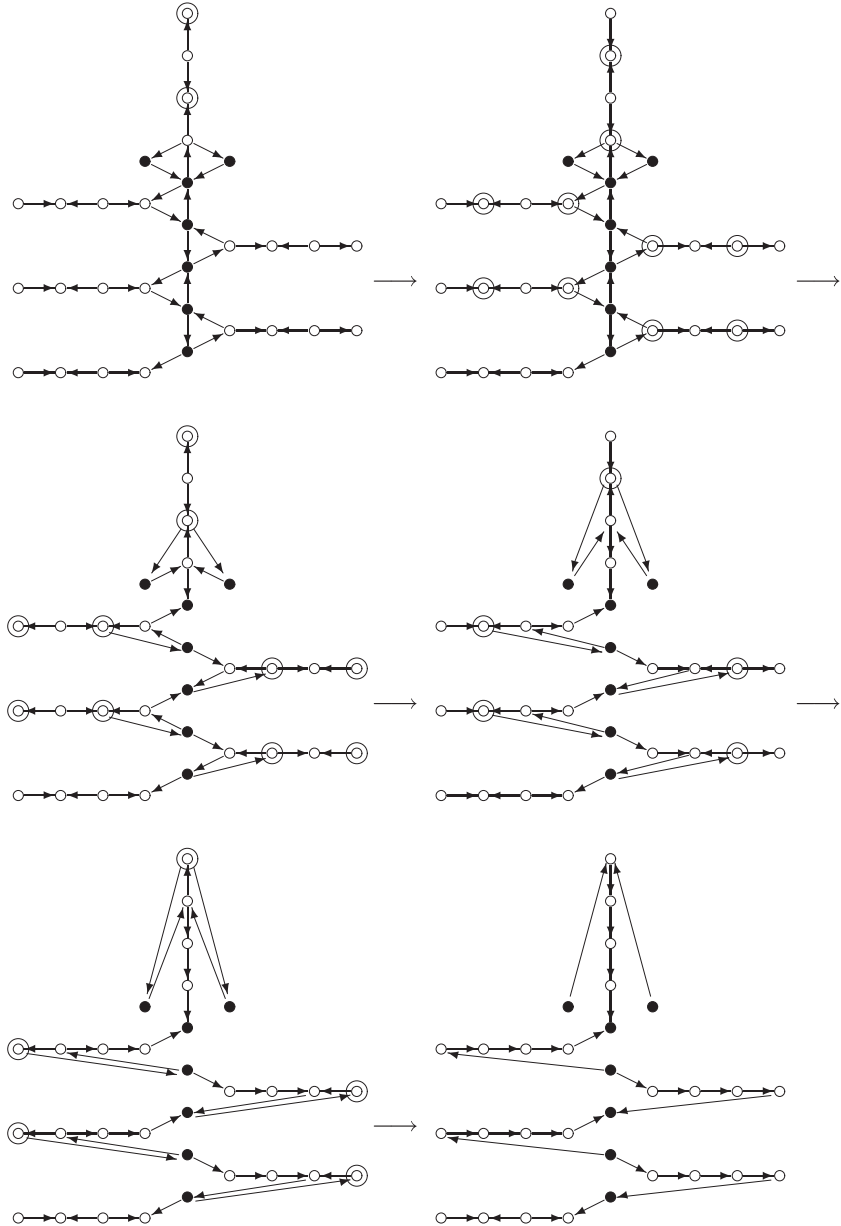}}
\caption{Sequence of
mutations from $Q_{4,7}(\mathrm{SG})$
to the quiver of type $D_{31}$.}
\label{fig:muteq}
\end{figure}

\begin{figure}[t]
\centerline{\includegraphics{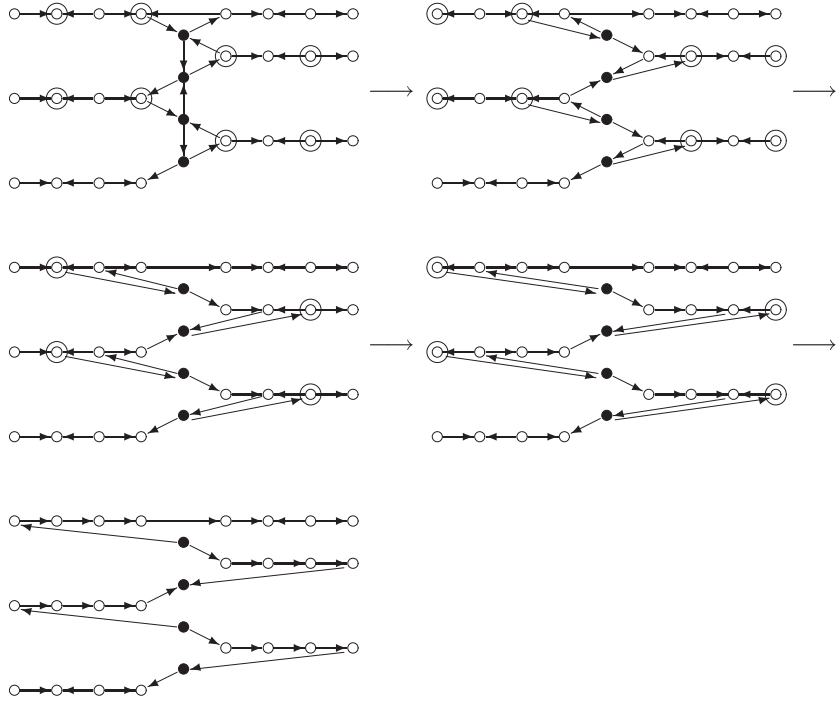}}

\caption{Sequence of
mutations from $Q_{4,7}(\mathrm{RSG})$
to the quiver of type $A_{28}$.}
\label{fig:muteqr}
\end{figure}

For example, for $(m,n)=(1,7)$, this is seen as follows:
$$
\centerline{\includegraphics{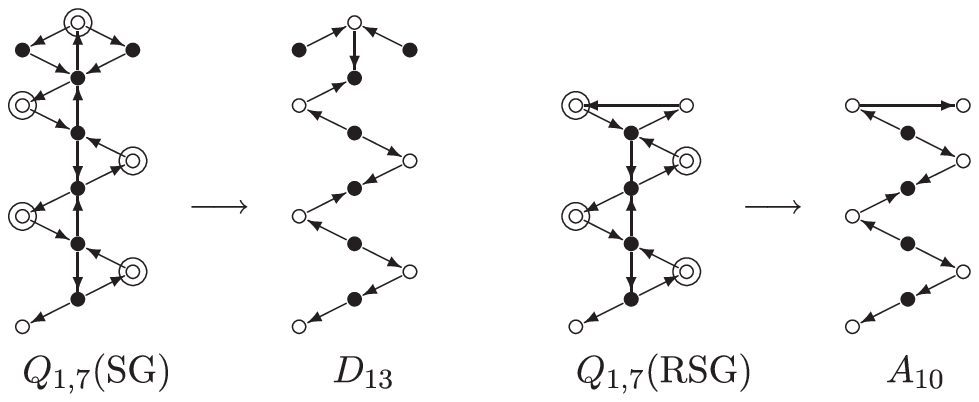}}
$$
Here, we do mutations at the encircled vertices.
For $m\geq 2$, the situation is a little more complicated
and the examples for $(m,n)=(4,7)$ are given in Figs.~\ref{fig:muteq} and~\ref{fig:muteqr}.

The referee made a further observation
on the periods in Theorems~\ref{thm:Yperiodm} and~\ref{thm:Yperiodmr}.

{\bf Fact 2.}
(a) Suppose that $mn-m+n$ is even. Then,
the period $2(mn-m+n)$ in Theo\-rems~\ref{thm:Yperiodm}$(a)$
 coincides with
 $h(D_{mn-m+n})+2$, which is the period
of the Coxeter mutation sequence for the cluster algebra
of type $D_{mn-m+n}$ studied by Fomin--Zelevinsky
in \cite{Fomin03a,Fomin03b,Fomin07}.
Similarly,
suppose that $mn-m+n$ is odd. Then,
the period $4(mn-m+n)$ in Theorems \ref{thm:Yperiodm}$(b)$
 coincides with
 $2(h(D_{mn-m+n})+2)$, which is the period
of the Coxeter mutation sequence for the cluster algebra
of type $D_{mn-m+n}$ in \cite{Fomin03a,Fomin03b,Fomin07}.

(b) The period $2(mn-m+n)$ in Theorems \ref{thm:Yperiodmr}
coincides with
 $2(h(A_{mn-m+n-3})+2)$, which is the period
of the Coxeter mutation sequence for the cluster algebra
of type $A_{mn-m+n-3}$ in~\cite{Fomin03a,Fomin03b,Fomin07}.

Facts 1 and 2 suggest that,
even though
the mutation sequences studied in this paper
and the ones studied in \cite{Fomin03a,Fomin03b,Fomin07}
are seemingly dif\/ferent,
they may be related, or may be treated in a unif\/ied way.
We leave this interesting question as an open problem.

\appendix
\section{Cluster algebras}

\label{sec:groc}

Here we collect basic def\/initions for cluster algebras
to f\/ix the convention and notation,
mainly following~\cite{Fomin07}.
For further necessary def\/initions and information
for cluster algebras, see~\cite{Fomin07}.

Let $I$ be a f\/inite index set throughout the appendix.

{\bf (i) Semif\/ield.}
A {\em semifield\/} $(\mathbb{P},\oplus, \cdot)$
is an Abelian multiplicative group
endowed with a binary
operation of addition $\oplus$ which is commutative,
associative, and distributive with respect to the
multiplication in $\mathbb{P}$.
The following three examples are important.

(a) {\em  Trivial semifield.}
The {\em trivial semifield\/}
$\mathbf{1}=\{\mathrm{1}\}$
consists of the multiplicative identity element $1$
with $1\oplus 1 = 1$.

\par
(b) {\em Universal semifield.}
For an $I$-tuple of variables $y=(y_i)_{i\in I}$,
the {\em universal semifield\/}
$\mathbb{P}_{\mathrm{univ}}(y)$
consists of  all the rational functions
of the form
$P(y)/Q(y)$  (subtraction-free  rational expressions),
where $P(y)$ and $Q(y)$ are the nonzero polynomials
 in $y_i$'s with {\em nonnegative\/} integer coef\/f\/icients.
The multiplication and the addition
are given by the usual ones of rational functions.

(c) {\em Tropical semifield.}
For an $I$-tuple of variables $y=(y_i)_{i\in I}$,
the {\em tropical semifield\/}
$\mathbb{P}_{\mathrm{trop}}(y)$
is the Abelian multiplicative group freely generated by
the variables  $y_i$'s endowed with the addition~$\oplus$
\begin{gather}
\label{eq:trop}
\prod_i y_i^{a_i}\oplus
\prod_i y_i^{b_i}
=
\prod_i y_i^{\min(a_i,b_i)}.
\end{gather}

{\bf (ii)  Mutations of matrix and quiver.}
An integer matrix
$B=(B_{ij})_{i,j\in I}$  is {\em skew-sym\-met\-ri\-zab\-le\/}
if there is a diagonal matrix $D=\mathrm{diag}
(d_i)_{i\in I}$ with $d_i\in \mathbb{N}$
such that $DB$ is skew-symmetric.
For a skew-symmetrizable matrix $B$ and
$k\in I$, another matrix $B'=\mu_k(B)$,
called the {\em mutation of $B$ at $k$\/}, is def\/ined by
\begin{gather}
\label{eq:Bmut}
B'_{ij}=
\begin{cases}
-B_{ij},& \mbox{$i=k$ or $j=k$},\\
B_{ij}+\frac{1}{2}
(|B_{ik}|B_{kj} + B_{ik}|B_{kj}|),
&\mbox{otherwise}.
\end{cases}
\end{gather}
The matrix $\mu_k(B)$ is also skew-symmetrizable.

It is standard to represent
a {\em skew-symmetric} (integer) matrix $B=(B_{ij})_{i,j\in I}$
by a {\em quiver $Q$
without loops or $2$-cycles}.
The set of the vertices of~$Q$ is given by~$I$,
and we put $B_{ij}$ arrows from~$i$ to~$j$
if $B_{ij}>0$.
The mutation $Q'=\mu_k(Q)$ of quiver $Q$ is given by the following
rule:
For each pair of an incoming arrow $i\rightarrow k$
and an outgoing arrow $k\rightarrow j$ in $Q$,
add a new arrow $i\rightarrow j$.
Then, remove a maximal set of pairwise disjoint 2-cycles.
Finally, reverse all arrows incident with $k$.

{\bf (iii)  Exchange relation of coef\/f\/icient tuple.}
Let $\mathbb{P}$ be a given semif\/ield.
For an $I$-tuple $y=(y_i)_{i\in I}$, $y_i\in \mathbb{P}$
and $k\in I$, another $I$-tuple $y'$ is def\/ined
by the {\em exchange relation}
\begin{gather}
\label{eq:coef}
y'_i =
\begin{cases}
\displaystyle
{y_k}{}^{-1}, &i=k,\\
\displaystyle
y_i \left(\frac{y_k}{1\oplus {y_k}}\right)^{B_{ki}},&
i\neq k,\ B_{ki}\geq 0,\\
y_i (1\oplus y_k)^{-B_{ki}},&
i\neq k,\ B_{ki}\leq 0.
\end{cases}
\end{gather}

{\bf (iv)   Exchange relation of cluster.}
Let $\mathbb{QP}$ be the  quotient f\/ield of the group ring
$\mathbb{Z}\mathbb{P}$ of $\mathbb{P}$,
 and let $\mathbb{QP}(z)$ be the rational function f\/ield of
algebraically independent variables, say, $z=(z_i)_{i\in I}$
over $\mathbb{QP}$.
For an $I$-tuple $x=(x_i)_{i\in I}$ which
is a free generating set of $\mathbb{QP}(z)$
and $k\in I$, another $I$-tuple $x'$ is def\/ined
by the {\em exchange relation}
\begin{gather}
\label{eq:clust}
x'_i =
\begin{cases}
{x_k}, &i\neq k,\\
\displaystyle
\frac{y_k
\prod\limits_{j: B_{jk}>0} x_j^{B_{jk}}
+
\prod\limits_{j: B_{jk}<0} x_j^{-B_{jk}}
}{(1\oplus y_k)x_k},
&
i= k.
\end{cases}
\end{gather}

{\bf (v)  Seed mutation.} For the above triplet $(B,x,y)$
in (ii)--(iv), which is called a {\em seed},  the mutation
$\mu_k(B,x,y)=(B',x',y')$  at $k$ is def\/ined
 by combining \eqref{eq:Bmut}, \eqref{eq:coef},
and \eqref{eq:clust}.

{\bf (vi)  Cluster algebra}. Fix a semif\/ield $\mathbb{P}$
and a seed ({\em initial seed\/}) $(B,x,y)$, where
$x=(x_i)_{i\in I}$  are algebraically independent variables
over $\mathbb{Q}\mathbb{P}$,
and
$y=(y_i)_{i\in I}$  are elements in $\mathbb{P}$.
Starting from $(B,x,y)$, iterate mutations and collect all the
seeds $(B',x',y')$.
We call  $y'$  and $y'_i$ a {\em coefficient tuple} and
a {\em coefficient}, respectively.
We call  $x'$  and $x'_i\in \mathbb{Q}\mathbb{P}(x)$, a {\em cluster} and
a {\em cluster variable}, respectively.
The {\em cluster algebra $\mathcal{A}(B,x,y)$ with
coefficients in $\mathbb{P}$} is a
$\mathbb{Z}\mathbb{P}$-subalgebra of the
rational function f\/ield $\mathbb{Q}\mathbb{P}(x)$
generated by all the cluster variables.

{\bf (vii) $\boldsymbol{F}$-polynomial.}
The cluster algebra $\mathcal{A}(B,x,y)$ with coef\/f\/icients in
the tropical semi\-f\/ield $\mathbb{P}_{\mathrm{trop}}(y)$ is called the
cluster algebra with the tropical coef\/f\/icients
(the principal coef\/f\/icients in~\cite{Fomin07}).
There, each cluster variable $x'_i$ is an element
in $\mathbb{Z}[x^{\pm 1},y]$.
The $F$-polynomial $F'_i(y)\in \mathbb{Z}[y]$ (for $x'_i$)
is def\/ined as the specialization of $x'_i$ with $x_i=1$ ($i\in I$).

\subsection*{Acknowledgements} It is our great pleasure to thank
the anonymous
referee who generously pointed out a crucial fact for the subject
of the paper.

\pdfbookmark[1]{References}{ref}
\LastPageEnding

\end{document}